%% file: SzymczakFunctorForRelations.tex
\documentclass[11pt]{amsart}
\usepackage{amsmath}
\usepackage{amscd}
\usepackage{pb-diagram}
\usepackage{comment}
\usepackage{graphicx}
\usepackage{amssymb}
\usepackage{pdfpages}
\usepackage{hyperref}

\input MM_MathDefs

\input MM_AlgDefs

\newcommand{\gdom}{\operatorname{gdom}}

\def\category#1{\text{\sc #1}}

\newcommand{\Szym}{\category{Szym}}
\newcommand{\FSet}{\category{FSet}}
\newcommand{\FRel}{\category{FRel}}
\newcommand{\LRel}{\category{LRel}}
\renewcommand{\Per}{\category{Per}}
\newcommand{\PerPts}{\operatorname{Per}}
\newcommand{\lab}{\operatorname{lab}}

\articletheorems

\begin{document}

\author{Mateusz Przybylski}
\address{Mateusz Przybylski, Division of Computational Mathematics,
  Faculty of Mathematics and Computer Science,
  Jagiellonian University, ul.~St. \L{}ojasiewicza 6, 30-348~Krak\'ow, Poland
}
\email{Mateusz.Przybylski@uj.edu.pl}
\author{Marian Mrozek}
\address{Marian Mrozek, Division of Computational Mathematics,
  Faculty of Mathematics and Computer Science,
  Jagiellonian University, ul.~St. \L{}ojasiewicza 6, 30-348~Krak\'ow, Poland
}
\email{Marian.Mrozek@uj.edu.pl}
\author{Jim Wiseman}
\address{Jim Wiseman, Department of Mathematics, Agnes Scott College, Decatur, GA 30030, USA
}
\email{jwiseman@agnesscott.edu}

\thanks{Research partially supported by
  the Polish National Science Center under Ma\-estro Grant No. 2014/14/A/ST1/00453
  and Opus Grant No. 2019/35/B/ST1/00874.}

\title[Szymczak Functor for Relations]{The Szymczak Functor on the Category of Finite Sets and Finite Relations}

\begin{abstract}
The Szymczak functor is a tool used to construct the Conley index for dynamical systems with discrete time. 
We present an algorithmizable classification of isomorphism classes in the Szymczak category over the category 
of finite sets with arbitrary relations as morphisms.  The research is the first step towards the construction of Conley theory for relations.
\end{abstract}

\maketitle


\section{Introduction}
\label{sec:intro}

In the 1970s Charles Conley proposed a homotopical invariant in dynamics \cite{Con1978},  
called after him the Conley index, which proved to be a very useful tool in the qualitative study of flows. 
The construction of the Conley index for flows is based on homotopies along the trajectories of the flow. 
This is an obstacle in the case of dynamical systems with discrete time, because in this case such homotopies do not make sense.
As a remedy several constructions based on shape theory and algebra were proposed \cite{RS1988,Mr1990,Mr1994} but it was
A.~Szymczak who indicated \cite{Sz1995} that all these constructions are functorial and are a special case of a certain universal functor associated with 
an arbitrary category $\cE$. To recall the functor consider the category $\Endo(\cE)$ of endomorphisms of $\cE$ whose objects are pairs $(E,e)$ 
with $e:E\to E$ an endomorphism in $\cE$ and whose morphisms are morphisms in $\cE$ commuting with the endomorphism.
Szymczak constructs a category $\Szym(\cE)$ and a functor $\Szym: \Endo(\cE)\to \Szym(\cE)$.
The feature of the Szymczak functor crucial in the construction of the Conley index is that it
assigns isomorphic objects in $\Szym(\cE)$  to objects $(A,a)$, $(B,b)$ in $\Endo(\cE)$ related by morphisms  
$\varphi:(A,a)\rightarrow (B,b)$ and $\psi:(B,b)\rightarrow (A,a)$ such that $b=\varphi\psi$, $a=\psi\varphi$.
Note that the Szymczak functor may be viewed as a functorial emanation of the concept of shift equivalence \cite{Wi1970},
another general concept used in the construction of the Conley index \cite{FR00}. Also, the Szymczak category can be interpreted as a localization of the $\Endo(\cE)$ category with respect to the class of morphisms $a\in \Endo(\cE)((A,a),(A,a))$ \cite{GZ1967}. 
The Szymczak functor is universal in the sense that any other functor with this feature factorizes through $\Szym$.
The universality of the Szymczak functor shows its generality but is also responsible for its computational weakness, 
because there is no general method to tell whether two objects in $\Szym(\cE)$ are isomorphic or not.  
In fact, it is even not clear when $\Szym(\cE)$ is not trivial, that is  whether $\Szym(\cE)$ contains non-isomorphic objects if $\cE$ does. 
For this reason, in concrete problems some less general but easier to compute functors are used, for instance the Leray functor \cite{RS1988,Mr1990,Mr1994}.
In rigorous algorithmic computations of the Conley index \cite{Sz1997,Mr2005} there is an additional challenge. 
Such computations, based on interval arithmetic \cite{Mo1966}, lead to multivalued dynamical systems and, in consequence, 
to categories whose morphisms are not maps but relations. 
Multivalued maps appear also in sampled dynamical systems constructed directly from data and acting on finite topological spaces \cite{DeyAtAl2019,DMS2020,BMW}. 
So far, the only method to deal with multivalued maps in the context of the Szymczak functor is to assume that
they have acyclic values,  because such maps induce single valued maps in homology. 
Acyclity may be achieved by enlarging the values at the expense of possible overestimation resulting in no interesting outcome.
However, the Szymczak category and the Szymczak functor are well defined for any category $\cE$, in particular for the category whose
morphisms are relations. Algorithmizable classification of isomorphism classes in Szymczak categories of relations 
could bring another method to study multivalued dynamical systems. In this paper we make a first step in this direction 
by providing such a method for the category of relations on finite sets.
The method, in particular, indicates that the Szymczak category for relations on finite sets is not trivial. 
Moreover, as we indicate in the last section, the proposed method proves non-triviality of Szymczak category for finite-dimensional vector spaces over a finite field with linear relations as morphisms. 

\section{Main results}
\label{sec:main}

Consider the category  $\FRel$ (see Section \ref{sec:szym_functor_in_frel}) consisting of finite sets as objects and binary relations as morphisms (arrows).
We interpret  an object $(X,R)$ in $\Endo(\FRel)$  as a directed graph with  $X$ as the set of vertices and $R$ as the set of edges.
We say that $(X,R)$ is  \emph{canonical} (see Subsection \ref{subsec:obj-in-canonical-form} for the precise definition)  if each vertex in $X$ belongs to a closed path, for each strongly connected
component $U\subset X$ the restriction $R_U:=R\cap U\times U$ is a bijection $R_U:U\to U$ and $R$ has periodic powers, that
is there exists a $p\geq 1$ such that $R^{p+1}=R$.

The following two theorems constitute the main theoretical results of the paper. 
We prove them in Section \ref{sec:szym_functor_in_frel}.

 \begin{thm}[see Theorem \ref{thm:each-obj-has-its-canonical-form}]
 \label{thm:main1}
   Every object in  $\Endo(\FRel)$ is isomorphic in  $\Szym(\FRel)$ to a canonical object.
 \end{thm}

 \begin{thm}[see Theorem \ref{thm:szym-isom-induces-can-obj-endo-isom}]
 \label{thm:main2}
Two canonical objects are isomorphic in  $\Szym(\FRel)$ if and only if they are isomorphic in $\Endo(\FRel)$.
 \end{thm}

\begin{table}
   \label{tab:class_summary}
  \begin{tabular}{ | c | c | c | }
    \hline
    $\card X$ & No. of objects & No. of $\Szym$ classes \\ \hline
    $\leq 1$ & $2$ & $2$ \\ 
    $\leq 2$ & $16$ & $5$ \\
    $\leq 3$ & $512$ & $14$ \\
    $\leq 4$ & $65536$ & $48$ \\
    $\leq 5$ & $33554432$ & $192$ \\
    \hline
   \end{tabular}   

   \caption{Number of different objects and different isomorphism classes  in $\Szym(\FRel)$
   for sets of cardinality not exceeding $n=1,2,3,4,5$.}   
\end{table}

\begin{figure}[h]
\begin{center}
  \includegraphics[width=0.5\textwidth]{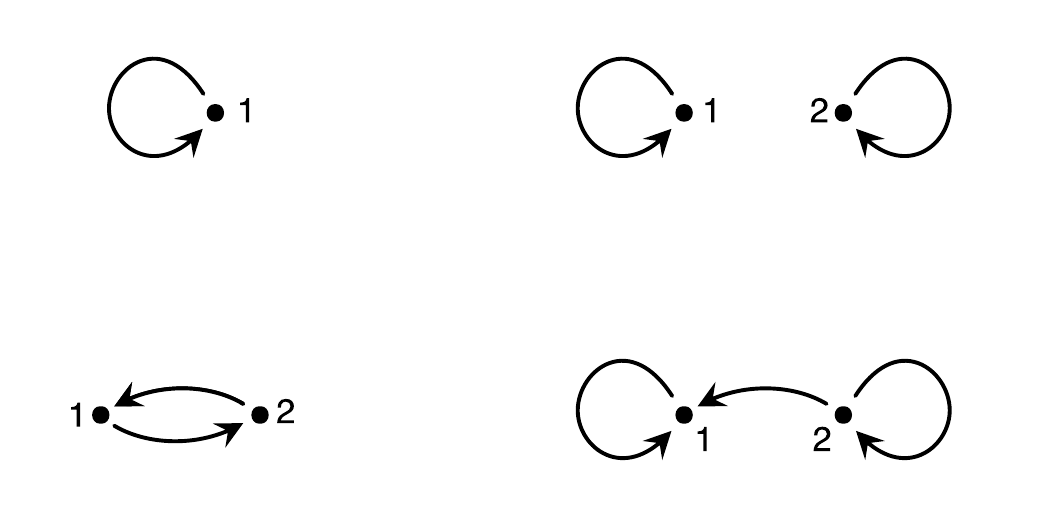}
\end{center}

  \caption{Canonical objects in $\Szym(\FRel)$ of cardinality zero (empty relation), one and two.
  }
  \label{fig:SzymClasses_012}
\end{figure}

\begin{figure}[h]
\begin{center}
  \includegraphics[width=0.95\textwidth]{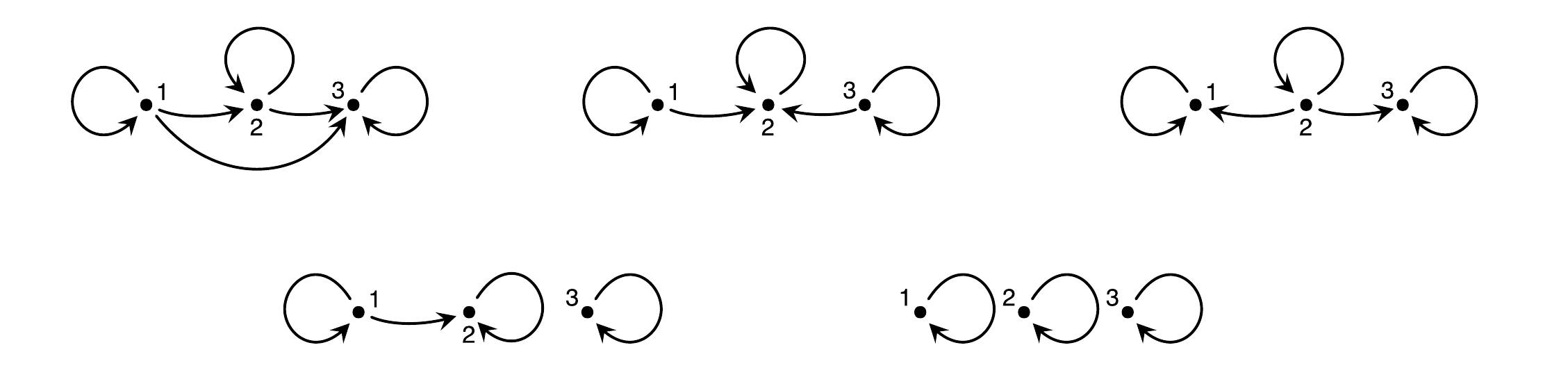}
\end{center}
  \caption{Canonical objects in $\Szym(\FRel)$ of cardinality three with three strongly connected components.
  }
  \label{fig:SzymClasses_3a}
\end{figure}

\begin{figure}[h]
\begin{center}
  \includegraphics[width=0.95\textwidth]{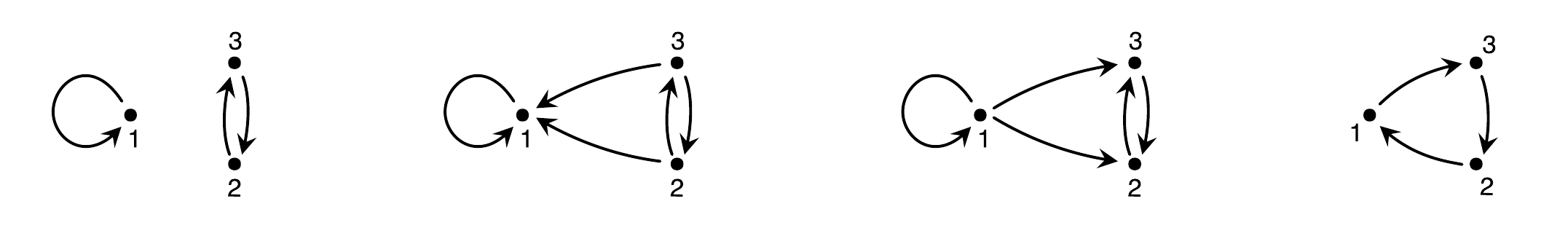}
\end{center}
  \caption{Canonical objects in $\Szym(\FRel)$ of cardinality three with less than three strongly connected components.
  }
  \label{fig:SzymClasses_3b}
\end{figure}

 Theorem~\ref{thm:main1} shows that each isomorphism class  in $\Szym(\FRel)$  admits a canonical representant. 
 Since the proof is constructive, the representant may be computed algorithmically. 
 Thus, the classification problem in $\Szym(\FRel)$ is reduced to the classical classification of graphs.
 This lets us compute canonical representants of isomorphism classes in $\Szym(\FRel)$ for sets of cardinality not exceeding five.  
 The number of different isomorphism classes is presented in Table~\ref{tab:class_summary}.
 The four canonical objects of cardinality zero one and two are presented in Figure~\ref{fig:SzymClasses_012}.
 The canonical objects of cardinality three are presented in Figures~\ref{fig:SzymClasses_3a}~and~\ref{fig:SzymClasses_3b}.

One can interpret relations on finite sets as Boolean matrices.  Then $(X,R)$ and $(Y,S)$ are isomorphic in $\Szym(\FRel)$ if and only if $R$ and $S$ are shift equivalent as Boolean matrices.  With some work, one can show that the linear algebraic result Proposition~3.5 from \cite{KR1986} (proven in \cite{KR1979}) is equivalent to part of Theorem~\ref{thm:main1} on canonical objects (the fact that any relation is isomorphic to a canonical form, though not the interpretation of that form).  The application in \cite{KR1979,KR1986} is to the classification of shifts of finite type, so there may be applications of Theorem~\ref{thm:main1} in that setting as well.


\section{Preliminaries}
\label{sec:pre}

We denote by $\NN_0$ ($\NN$, for short)  and $\NN_1$ the set of natural numbers respectively including and excluding zero.
We recall that a {\em binary relation} in $X\times Y$, or briefly a {\em relation},
is a subset $R\subset X\times Y$.
For a relation $R$ we use the classical notation $xRy$ to denote $(x,y)\in R$.
If $X'\subset X$ and $Y'\subset Y$, we call the relation $R_{|X'\times Y'}:=R\,\cap\, X'\times Y'$
the {\em restriction} of $R$ to $X'\times Y'$.
If $X=Y$ we say that $R$ is a relation in $X$.
If  $A\subset X$, by the restriction of $R$ to $A$
we mean the restriction of $R$ to $A\times A$.
We denote this restriction by $R_{|A}:=R\,\cap\, A\times A$.

The \emph{inverse relation} of $R$ is
\[
   R^{-1}:=\setof{(y,x)\in Y\times X\mid xRy}.
\]
Given another relation $S\subset Y\times Z$  we define the \emph{composition} of $S$ with $R$ as the relation
\[
   S\circ R:=\setof{(x,z)\in X\times Z\mid\text{ $xRy$ and $yRz$ for some $y\in Y$}}.
\]

The following proposition follows immediately from the definition of composition of relations.
\begin{prop}
\label{prop:rel-comp-monotonicity}
If $S\subset R\subset X\times X'$ and $S'\subset R'\subset X'\times X''$, then
$S'\circ S\subset R'\circ R$.
\end{prop}

The \emph{domain} of $R$ is
\[
  \dom R:=\setof{x\in X\mid \text{ $xRy$ for some $y\in X$}}
\]
and the \emph{image} of $R$  is
\[
  \im R:=\setof{y\in Y\mid \text{ $xRy$ for some $x\in X$}}.
\]

The following proposition is straightforward.
\begin{prop}
\label{prop:dom-im-comp}
Let $R\subset X\times Y$ and $S\subset Y\times Z$ be relations. Then
\[
    \dom S\circ R\subset \dom R \text{ and }  \im S\circ R\subset \im S.
\]
\qed
\end{prop}

The \emph{identity relation} on $X$ is
\[
  \id_X=\setof{(x,x)\in X\times X\mid x\in X}.
\]
For $n\in\ZZ$ the \emph{$n$th power} of a relation $R$ in $X$ is given recursively by
\[
   R^n:=\begin{cases}
           \id_X & \text{ for $n=0$,}\\
           R\circ R^{n-1} & \text{ for $n>0$,}\\
           R^{-1}\circ R^{n+1} & \text{ for $n<0$.}\\
        \end{cases}
\]

For $A\subset X$ we define the \emph{image of $A$ in $R$} as
\[
  R(A):=\setof{y\in Y\mid\text{ $xRy$ for some $x\in A$}}.
\]
In case $A$ is a singleton we simplify notation and write $R(x)$ meaning $R(\{x\})$.
In particular $xRy$ if and only if $y\in R(x)$.

With every relation $R\subset X\times Y$ we can associate the map
\[
  X\ni x \mapsto R(x)\in 2^X.
\]
We prefer to think of the value $R(x)$ not as an element of the space of subsets $2^X$
but as a subset of $X$, i.e. we consider $R$ as a \emph{partial multivalued map}
\[
  X\supset \dom R\ni x \mapsto R(x)\subset X.
\]
The change of terminology serves to emphasize that we want to think of $R$
as a map which may have many values.
In particular, in the sequel we will write $y\in R(x)$ instead of $xRy$
to emphasize the multivalued map interpretation of the relation $R$.
Note that the value $R(x)$ may be empty. This happens when $x\not\in\dom R$.
We say that a relation $R$ is a \emph{multivalued map} and we write $R: X\mto Y$ if $\dom R=X$.
If $R(x)$ is a singleton for every $x\in\dom R$ we identify $R(x)$ with its unique value
and we say that the relation $R$ is a \emph{partial map} and we write $R: X\pto Y$. 
We say that a partial map $R\subset X\times Y$ is a \emph{map} and we write $R: X\to Y$ 
if additionally $\dom R=X$.

We say that a relation $R\subset X\times Y$ is \emph{injective} if $R(x_1)\cap R(x_2)\neq\emptyset$
implies $x_1=x_2$ for any $x_1,x_2\in\dom R$. We say that a relation $R\subset X\times Y$ is \emph{surjective}
if $\im R=Y$.
We say that $g\subset X\times Y$ is a \emph{bijection} or a \emph{bijective map} if it is an injective and surjective map.
Note that  a relation which is both injective and surjective need not be
a bijection or even a map. But, we have the following proposition.
\begin{prop}
\label{prop:bijective-composition}
Let $R\subset X\times Y$ be a relation and let $S\subset Y\times Z$ be a multivalued map, that is $\dom S=Y$.
If $S\circ R\subset X\times Z$
is a bijective map then $S$ is a surjective  multivalued map and $R$ is an injective multivalued map.
\end{prop}
\proof
  Let $g:=S\circ R$.  Since $g$ is a bijection, we have  $\im g = Z$ and $\dom g= X$.
  It follows from Proposition~\ref{prop:dom-im-comp} that $Z=\im g=\im S\circ R\subset \im S$.
  Hence, $\im S=Z$ which means that $S$ is a surjection.
  Similarly,  $X=\dom g=\dom S\circ R\subset \dom R$.
  Hence, $\dom R=X$ which means that $R$ is a multivalued map.
  To see that $R$ is injective assume that $R(x_1)\cap R(x_2)\neq\emptyset$.
  Let $y\in R(x_1)\cap R(x_2)$.  Since $\dom S=Y$, we can find a $z\in Z$ such that
  $ySz$. It follows that $x_1\, S\circ R\, z$ and $x_2\, S\circ R\, z$.
  Since $g=S\circ R$ is a bijection we obtain $x_1=x_2$.
\qed

By a {\em directed graph} (or just a {\em digraph}) we mean a pair $G:=(V,E)$ consisting of the finite set of vertices $V$ and the set of edges $E\subset V\times V$. We allow a digraph to contain loops, that is edges in the form of $(v,v)$, where $v\in V$. A {\em walk} in $G$ is a sequence $x=x_1x_2\ldots x_k$ such that $x_i\in V$ for $i=1,\ldots,k$ and $(x_i,x_{i+1})\in E$ for $i=1,\ldots,k-1$. We then say that $x$ is a {\em walk from $x_1$ to $x_k$} or just an {\em $(x_1,x_k)$-walk}. The {\em length} of walk $x$ is the number of edges $(x_i,x_{i+1})$, that is $k-1$. We denote it by $\#x$. We say that a vertex $x_i$ {\em lies on a walk $x$} if it is contained in the sequence that constitutes the walk $x$. If the vertices of the walk $x$ are different, then we call $x$ a {\em path} (or a {\em path from $x_1$ to $x_k$}). A walk $x=x_1\ldots x_k$ of length greater or equal to $1$ is a {\em cycle} if $x_1 = x_k$. A {\em concatenation of a walk $x=x_1\ldots x_k$ with a walk $y=y_1\ldots y_n$} is a walk $xy:=x_1\ldots x_ky_2\ldots y_n$ provided $x_k=y_1$.  

A digraph $G=(V,E)$ is {\em strongly connected} if for each $v,u\in V$ there exist a $(v,u)$-walk and a $(u,v)$-walk and both walks have length greater than $0$. For any digraph $G=(V,E)$ a set $U\subset V$ is called a {\em strongly connected component of $G$} if the digraph $G(U):=(U,\{(v,u)\in E\ |\ v,u\in U\})$ is strongly connected and there is no other $W$ such that $U\subset W\subset V$ and $G(W)$ is strongly connected. In this paper we do not use any other connectivity of digraphs. Therefore, we often use a brief form connected meaning strongly connected.

Each relation $R\subset X\times X$ may be considered as the directed graph $(X,R)$. Similarly, any directed graph $G=(V,E)$ may be considered as the binary relation $E\subset V\times V$. This observation lets us use the notions of digraph and relation interchangeably throughout the paper, choosing the one that better fits to the presented content and applying digraph terminology to relations and vice versa. 

Notice that the existence of a $(c_1,c_p)$-walk of length $p-1$ in the digraph $(X,R)$ is equivalent to the fact $c_p\in R^{p-1}(c_1)$. In particular, if $(c_1,c_p)$-walk is a cycle, then the existence of a cycle is equivalent to $c_i\in R^{p-1}(c_i)$ for each $i=1,\ldots,p$.

\section{Szymczak functor}

\subsection{Categories.}
Let $\cE$ be a category. Recall that a morphism $\varphi:E\to E'$ is an isomorphism
in $\cE$ if there exists a morphism $\psi:E'\to E$ such that $\psi\circ \varphi=\id_E$
and $\varphi\circ\psi=\id_{E'}$. Then also $\psi$ is an isomorphism. It is uniquely
determined by $\varphi$ and called the inverse morphism of $\varphi$. We denote it $\varphi^{-1}$.
We recall that an \emph{endomorphism}   in $\cE$ is a morphism of the form $e:E\to E$,
that is a morphism whose source object is the same as the target object.
An \emph{automorphism} is an endomorphism which is also an isomorphism.

Let $E$ and $F$ be two objects of $\cE$ and let $e:E\to E$, $f:F\to F$
be morphisms in $\cE$.
We say that $e$ and $f$
are {\em conjugate} if there exists an isomorphism $\varphi:E\to F$ such that
$\varphi\circ e=f\circ \varphi$.

\begin{prop}
\label{prop:four-iso}
Assume  the diagram
\[
\begin{diagram}
 \node{E}
 \arrow{s,l}{\varphi}
 \arrow{e,t}{\mu}
 \node{E}
 \arrow{s,r}{\varphi}
 \\
 \node{F}
 \arrow{e,t}{\nu}
 \arrow{ne,t}{\psi}
 \node{F}
\end{diagram}
\]
of morphisms in $E$ is commutative. If $\mu$ and $\nu$ are isomorphisms, then so are $\varphi$ and $\mu$.
In particular, the isomorphisms $\mu$ and $\nu$ are conjugate.
\end{prop}
\proof
  Set $\varphi':=\varphi\circ\mu^{-1}$. Then $\psi\circ\varphi'=\id_E$.
From $\nu\circ\varphi=\varphi\circ\mu$ we get $\varphi\circ\mu^{-1}=\nu^{-1}\circ\varphi$.
Therefore,   $\varphi'\circ\psi=\varphi\circ\mu^{-1}\circ\psi=\nu^{-1}\circ\varphi\circ\psi=\nu^{-1}\circ\nu=\id_F$.
This proves that $\psi$ is an isomorphism. It follows that $\varphi=\psi^{-1}\circ\mu$ is an isomorphism
as a composition of isomorphisms.
\qed

\subsection{Category of endomorphisms.} We define the category of endomorphisms of $\cE$,
denoted by $\Endo(\cE)$, as follows:  the objects of
$\Endo(\cE)$ are pairs $(E,e)$, where $E\in\cE$ and $e\in\cE(E,E)$
is an endomorphism of $E$.  The set of morphisms from
$(E,e)\in \Endo(\cE)$ to $(F,f)\in \Endo(\cE)$ is the subset of $\cE(E,F)$
consisting of exactly those morphisms $\varphi\in\cE(E,F)$ for which $f\varphi=\varphi
e$.
We write $\varphi:(E,e)\rightarrow (F,f)$ to denote that $\varphi$ is a morphism from $(E,e)$
to
$(F,f)$ in $\Endo(\cE)$.  Note that in particular $e:(E,e)\map (E,e)$ is e morphism in   $\Endo(\cE)$.

Let $\cC$ be another category and let $L:\Endo(\cE)\map\cC$ be a functor.
We say that $L$ is {\em normal} if $L(e)$ is an isomorphism in $\cC$
for any endomorphism $e:E\map E$ in $\cE$. We have the following theorem.

\begin{thm}
\label{thm:normal-functor}
Assume $L:\Endo(\cE)\map\cC$ is a normal functor and
$\varphi:(E,e)\rightarrow (F,f)$, $\psi:(F,f)\rightarrow (E,e)$
are such that $e=\varphi\psi$, $f=\psi\varphi$.  Then we have the commutative diagram
\[
\begin{diagram}
 \node{L(E,e)}
 \arrow{s,l}{L(\varphi)}
 \arrow{e,t}{L(e)}
 \node{L(E,e)}
 \arrow{s,r}{L(\varphi)}
 \\
 \node{L(F,f)}
 \arrow{e,t}{L(f)}
 \arrow{ne,t}{L(\psi)}
 \node{L(F,f)}
\end{diagram}
\]
in $\cC$, in which all morphisms are isomorphisms.
\end{thm}
\proof
The theorem is an immediate consequence of Proposition~\ref{prop:four-iso}.
\qed

\subsection{Szymczak category.} With every category $\cE$ one can associate its Szymczak category $\Szym(\cE)$
defined as follows. The objects of $\Szym(\cE)$ are the objects of $\Endo(\cE)$.
Given objects $(E,e)$ and $(E',e')$ in $\Szym(\cE)$
we consider the equivalence relation
in $\Endo(\cE)((E,e),(E',e'))\times\NN$ defined by
\[
(\varphi,m)\equiv (\varphi',m')
\]
for $(\varphi,m),(\varphi',m')\in \Endo(\cE)((E,e),(E',e'))\times\NN$
if and only if there exists a $k\in\NN$ such that
\begin{equation}
\label{szym-equiv}
   \varphi \circ e^{m'+k} = \varphi' \circ e^{m+k}.
\end{equation}
We define the set of morphisms
$\Szym(\cE)((E,e),(E',e'))$ as the collection of equivalence classes of the relation $\equiv$.
Given morphisms $[\varphi,m]:(E,e)\to(E',e')$ and $[\varphi,m']:(E',e')\to(E'',e'')$ we define
their composition by
\[
    [\varphi',m']\circ[\varphi,m]:=[\varphi'\circ\varphi,m+m']
\]
One easily verifies that the composition is well defined and $[\id_E,0]$
is the identity morphism on $(E,e)$. Thus, $\Szym(\cE)$ is indeed a category.

   There is a functor $\Szym:\Endo(\cE)\map\Szym(\cE)$ which fixes objects
and sends a morphism $\varphi:(E,e)\map (F,f)$ to the equivalence class $[\varphi,0]$.
In general, it may happen that $\Szym(\varphi)=\Szym(\varphi')$ even if $\varphi\neq\varphi'$.
Nevertheless it is convenient to write just $\varphi$ to denote $\Szym(\varphi)$ whenever
it is clear from the context in which category we work.
One easily verifies that every morphism  $e:(E,e)\map(E,e)$ in $\Szym(\cE)$ has an inverse given by
\[
   \bar{e}:=[\id_E,1].
\]
Indeed, we have
\[
   e\circ\bar{e}=[e,0]\circ[\id_E,1]=[e,1]=[\id_E,0]=\id_{(E,e)}
\]
which shows that $\bar{e}$ is an inverse of $e$. We can also write the abstract morphism $[\varphi,n]$
in terms of $\bar{e}$ as
\begin{equation}
\label{eq:szym-inverse}
     [\varphi,n]=[\varphi,0]\circ[\id_E,1]^n=\varphi\circ\bar{e}^n.
\end{equation}

Thus,  $\Szym(e)$ is invertible in $\Szym(\cE)$.  Therefore, $\Szym$ is a normal functor.
Actually, this is the most general normal functor in the following sense.
\begin{thm}
\label{thm:Szym-universal}  \cite[Theorem 6.1]{Sz1995}
For every normal functor $L:\Endo(\cE)\rightarrow \cC$ there exists
 a unique functor $L':\Szym(\cE)\rightarrow \cC$ such that the diagram
\[
\begin{diagram}
 \node{\Endo(\cE)}
 \arrow{s,l}{\Szym}
 \arrow{e,t}{L}
 \node{\cC}
 \\
 \node{\Szym(\cE)}
  \arrow{ne,b}{L'}
\end{diagram}
\]
commutes.
\end{thm}

The construction of the Szymczak category and the Szymczak functor
is due to Szymczak \cite{Sz1995}.

We say that two objects  $(E,e)$ and $(E',e')$ of $\Endo(\cE)$ are \emph{conjugate}  if $e$ and $e'$ are conjugate in $\cE$.
\begin{prop}
\label{prop:Szym-conj}
Assume $(E,e)$ and $(E',e')$ are conjugate objects of $\Endo(\cE)$. 
Then  $(E,e)$ and $(E',e')$ are isomorphic in  $\Szym(\cE)$.
\end{prop}
\proof
   Let $\varphi:E\to E'$ be an isomorphism in $\cE$ such that $\varphi\circ e=e'\circ \varphi$ and let $\psi:=\varphi^{-1}$.
   Then $[\psi,0]\circ[\varphi,0]=[\id_E,0]$ and $[\varphi,0]\circ[\psi,0]=[\id_{E'},0]$,
   which proves that $(E,e)$ and $(E',e')$ are isomorphic in  $\Szym(\cE)$.
\qed

It is not difficult to give examples showing that the converse of Proposition~\ref{prop:Szym-conj} is not true.
However, it is true in the category $\Auto(\cE)$ defined as the full subcategory of $\Endo(\cE)$ whose objects
are objects $(E,e)$ of $\Endo(\cE)$ such that $e$ is an isomorphism in $\cE$.
Indeed, we have the following proposition.
\begin{prop}
\label{prop:Auto-conj}
Assume $(E,e)$ and $(E',e')$ are objects in $\Auto(\cE)$.
If $\Szym(E,e)  \cong \Szym(E',e')$, then $(E,e)$ and $(E',e')$ are conjugate.
\end{prop}
\proof
   Since $\Szym(E,e)  \cong \Szym(E',e')$, we may find morphisms $\varphi:(E,e)\to (E',e')$
   and $\psi:(E',e')\to (E,e)$ as well as constants $n,n'\in\NN_0$ such that
   $[\varphi,n]\circ[\psi,n']=[\id_E,0]$ and $[\psi,n']\circ[\varphi,n]=[\id_{E'},0]$.
   This means that there exist $k,k'\in\NN_0$ such that
   $\psi\circ\varphi\circ e^k=e^{k+n+n'}$ and $\varphi\circ\psi\circ e'^{k'}=e'^{k'+n+n'}$.
   Since $e$ and $e'$ are isomorphisms, the equalities may be reduced to
   $\psi\circ\varphi=e^{n+n'}$ and $\varphi\circ\psi=e'^{n+n'}$.
   Since both $e^{n+n'}$ and $e'^{n+n'}$ are isomorphisms, the conclusion follows
   now immediately from Proposition~\ref{prop:four-iso}.
\qed

The Szymczak category can be seen as a localization of the $\Endo(\cE)$ category with respect to the class of morphisms $e\in\Endo(\cE)((E,e),(E,e))$ (see \cite{GZ1967}).

The Szymczak category and the Szymczak functor are very general concepts, defined for any category.
However, in practical terms it is not obvious how to compute the Szymczak category and Szymczak functor for concrete categories.
In the next section we do it for the category of finite sets.


\section{Szymczak functor in $\FSet$}

Let $\FSet$ denote the category of finite sets with maps as morphisms.
Given an object $(X,f)$ in $\Endo(\FSet)$ we say that an $x\in X$ is a {\em periodic point} of $f$
if there exists a $k\in\NN_1$ such that $f^n(k)=x$.
We then say that $k$ is a {\em period} of $x$ and $x$ is $k$-periodic.
We denote the set of periodic points of $f$ by $\PerPts f$
and the set of $k$-periodic points of $f$ by $\PerPts^k f$.

Let $(X,f)$ be a fixed object of $\Endo(\FSet)$.

\begin{prop}
\label{prop:f-Per}
The map
\[
   f_{|\PerPts f}:\PerPts f\ni x\mapsto f(x)\in\PerPts f
\]
is a well defined bijection.
\end{prop}
\proof
  Obviously, if $x$ is $k$-periodic than so is $f(x)$.
Therefore, the map is well defined.
Assume $x\in\PerPts f$ and let $k$ be a period of $x$.
Then $f^{k-1}(x)\in\PerPts f$ and $x=f(f^{k-1}(x))$, which proves that $f$ is a surjection.
To see that $f$ is an injection take $x_1,x_2\in\PerPts f$ such that $f(x_1)=f(x_2)$.
Let $k_i$ be a period of $x_i$ for $i=1,2$ and let $k:=k_1k_2$.
Then $x_1=f^k(x_1)=f^{k-1}(f(x_1))=f^{k-1}(f(x_2))=f^k(x_2)=x_2$.
\qed

Consider the relation $\sim_f$ in $X$ defined by
\begin{equation}
\label{eq:f-equiv}
  x\sim_fy \text{ if and only if } \exists n\in\NN_0\\:f^n(x)=f^n(y).
\end{equation}
One easily verifies that $\sim_f$ is an equivalence relation in $X$.
Denote by $[x]_f$ the equivalence class of $x$ in $\sim_f$.

By a {\em periodic exponent} of $f$ we mean any multiplicity of  minimal periods of all points in $\PerPts f$
which is not less than the cardinality of $X$.

\begin{thm}
\label{thm:fn-Per-f}
Let $n$ be a periodic exponent of $f$. Then
$[x]_f\cap\PerPts f=\{f^n(x)\}$ for every $x\in X$.
\end{thm}
\proof
Denote by $g:=f_{|\PerPts f}^{-1}:\PerPts f\to\PerPts f$ the inverse of $f_{|\PerPts f}$
in $\PerPts f$. The inverse exists by Proposition~\ref{prop:f-Per}.
Fix an $x\in X$.

Since $X$ is finite, there exist a $k\in\NN_1$ and an $m\in\NN_0$ such that
  $f^{k+m}(x)=f^m(x)$. It follows that the set
\[
   J_x:=\setof{i\in\NN_0\mid f^i(x)\in\PerPts f}
\]
is non-empty.
Let $k_x:=\min J_x$ and let $u_x:=g^{k_x}(f^{k_x}(x))$.
Obviously $u_x\in\PerPts f$.
Since $f^{k_x}(u_x)=f^{k_x}(g^{k_x}(f^{k_x}(x)))=f^{k_x}(x)$, we see that $u_x\in [x]_f$.
Therefore $u_x\in [x]_f\cap\PerPts f$. We will prove that
\begin{equation}
\label{eq:fn-Per-f}
u_x=f^n(x).
\end{equation}
Since $n$ is a multiple of the minimal period of $u_x$, we have $f^n(u_x)=u_x$.
Note that $n\geq k_x$. 
Therefore,
\begin{eqnarray*}
u_x&=&f^n(u_x)=f^n(g^{k_x}(f^{k_x}(x)))=f^{n-k_x}(f^{k_x}(g^{k_x}(f^{k_x}(x))))=\\
   &=&f^{n-k_x}(f^{k_x}(x))=f^n(x),
\end{eqnarray*}
which proves \eqref{eq:fn-Per-f}. It follows that $\{f^n(x)\}\subset [x]_f\cap\PerPts f$.
To prove that the two sets are actually equal, it suffices to show that
$[x]_f\cap\PerPts f$ contains at most one point.
For this end assume that $x_1,x_2\in [x]_f\cap\PerPts f$.
Then there exist $n_1,n_2\in\NN_0$ such that $f^{n_1}(x_1)=f^{n_1}(x)$
and $f^{n_2}(x_2)=f^{n_2}(x)$.
Since $x_1,x_2\in\PerPts f$, there exist $k_1,k_2\in\NN_1$ such that $f^{k_1}(x_1)=x_1$ and  $f^{k_2}(x_2)=x_2$.
Choose a $p\in\NN_1$ such that $m:=pk_1k_2\geq\max(n_1,n_2)$.
It follows that
\[
   x_1=f^m(x_1)=f^m(x)=f^m(x_2)=x_2.
\]
\qed

\begin{cor}
\label{cor:fn-Per-f}
Let $n$ be the minimal periodic exponent of $f$.
Then $f^m(X)= \PerPts f$ for $m\geq n$.
\end{cor}
\proof
  It follows from Theorem~ \ref{thm:fn-Per-f} that
  $f^n(X)\subset \PerPts f$. Therefore,   $f^m(X)\subset f^n(X)\subset \PerPts f$  for $m\geq n$.
  The opposite inculsion is obvious, becasue a periodic point of $f$ is in the image of $f^m$ for
  every $m\in\NN_0$.
\qed

\begin{cor}
\label{cor:per-exp}
Let $n_1$ and $n_2$ be periodic exponents of $f$.
Then $f^{n_1}=f^{n_2}$.
\end{cor}
\proof
Without loss of generality we may assume that $n_1$ is the minimal periodic exponent of $f$.
Let $k:=n_2-n_1$ and let $x\in X$. It follows from Corollary~\ref{cor:fn-Per-f} that
$f^{n_1}(x)\in\PerPts f$. Since $k$, as a difference of periodic exponents, is a multiple of
minimal periods of all points in  $\PerPts f$, we get
$f^{n_2}(x)=f^k(f^{n_1}(x))=f^{n_1}(x)$.
\qed

We refer to the common value of $f^n$ with $n$ being a periodic exponent of $f$
as the {\em universal power} of $f$.  We denote it $\hat{f}$.
    Note that by Corollary~\ref{cor:fn-Per-f}  we can consider $\hat{f}$ as a map $\hat{f}:X\to\Per f$.

\begin{prop}
\label{prop:universal-power}
For every $x\in\PerPts f$ we have $\hat{f}^{-1}(x)=[x]_f$.
In particular, the map
\[
  \hat{f}:X\ni x\mapsto \hat{f}(x)\in\PerPts f
\]
is a well defined surjection.
\end{prop}
\proof
Let $x'\in[x]_f$.  Then $f^k(x)=f^k(x')$ for some $k\in\NN_0$.
Let $n$ be a periodic exponent of $f$ which is greater then $k$.
Then $x=f^n(x)=f^n(x')=\hat{f}(x')$ and $x'\in \hat{f}^{-1}(x)$.
This proves that $[x]_f\subset \hat{f}^{-1}(x)$. To prove the opposite inclusion
take an $x'\in\hat{f}^{-1}(x)$. Then $\hat{f}(x')=x$.
Let $n$ be a periodic exponent of $f$.
It follows that $f^n(x')=\hat{f}(x')=x=f^n(x)$.
In consequence, $x'\in[x]_f$.
\qed

\begin{prop}
\label{prop:universal-representant}
Assume $(X,f)$ is an object of $\Endo(\FSet)$.
Let $\iota:\PerPts f\to X$ denote the inclusion map
and let $n$ be a periodic exponent of $f$.
Then, 
\[
[\iota,0]:(\PerPts f,f_{|\PerPts f})\to (X,f)
\]
and 
\[
[\hat{f},n]:(X,f)\to (\PerPts f,f_{|\PerPts f})
\]
are mutually inverse isomorphisms in $\Szym(\FSet)$. 
Moreover,
the map $\hat{f}$ induces a bijection $\bar{f}:(X/\sim_f,f')\to(\PerPts f,f_{|\PerPts f})$.
\end{prop}
\proof
   The equality   $\hat{f}=f^n=\id_X\circ f^n$ implies that 
\[
   [\iota,0]\circ [\hat{f},n]= [\hat{f},n]=[\id_X,0].
\]
Since  $f(\Per f)\subset\Per f$, we also have 
\[
\hat{f}_{|\Per f}=(f^n)_{|\Per f}=f_{|\Per f}^n=\id_{|\Per f}\circ f_{|\Per f}^n,
\]
which  implies  
\[
   [\hat{f}_{|\Per f},n]=[\id_{|\Per f},0].
\]
This proves that  $[\iota,0]$ and $[\hat{f},n]$
are mutually inverse isomorphisms in $\Szym(\FSet)$. 
It follows from the definition \eqref{eq:f-equiv} of the equivalence relation $\sim_f$ that $\bar{f}$ is well defined.
Moreover,  Corollary~\ref{cor:fn-Per-f} implies that $\hat{f}$ is a surjection. Hence, so is $\bar{f}$.
In consequence,   $\bar{f}$ is a bijection, because $X$ is finite.
\qed

Proposition~\ref{prop:f-Per} lets us define a functor
\[
\Per:\Endo(\FSet)\to\Auto(\FSet)
\]
as follows. For an object $(X,f)$ in $\Endo(\FSet)$ we set $\Per(X,f):=(\PerPts f,f_{|\PerPts f})$.
Given a morphism $\varphi:(X,f)\to(X',f')$
 we define $\Per(\varphi)$
as the map $\PerPts f\ni x\mapsto \varphi(x)\in\PerPts f'$. Note that this map is
well defined, because $x\in \PerPts^kf$ implies $f'^{k}(\varphi(x))=\varphi(f^k(x))=\varphi(x)$.
One easily verifies that $\Per$ is indeed a functor.
Moreover, it is a normal functor, because $\Per(f)$, as a bijection, is an isomorphism
in $\Auto(\FSet)$.

Let $\Per':\Szym(\FSet)\to \Auto(\FSet)$ be the functor associated to $\Per$ by Theorem~ \ref{thm:Szym-universal}. In particular, we have
\begin{equation}
\label{eq:PerSzym}
\Per'\circ\Szym=\Per.
\end{equation}

\begin{thm}
\label{thm:f-Per}
The functor $\Per'$ is a bijector.
\end{thm}
\proof
  We need to show that $\Per'$ is an injector and a surjector.
For this end assume $[\varphi,n]:\Szym(X,f)\to\Szym(X',f')$ and $[\psi,m]:\Szym(X,f)\to\Szym(X',f')$
are morphisms in $\Szym(\FSet)$ such that
\[
   \Per'([\varphi,n]) =\Per'([\psi,m]).
\]
Rewriting this formula using the functoriality of $\Per'$,
\eqref{eq:szym-inverse}, \eqref{eq:PerSzym}   and multiplying on the right
by $\Per'(f)^{m+n}$  we obtain
\begin{eqnarray}
\Per'(\varphi\circ\bar{f}^n) &=&\Per'(\psi\circ\bar{f}^m),\\
\Per'(\varphi)\circ\Per'(\bar{f})^n &=&\Per'(\psi)\circ\Per'(\bar{f})^m,\\
\Per'(\varphi)\circ\Per'(f)^m &=&\Per'(\psi)\circ\Per'(f)^n,\\
\Per(\varphi)\circ\Per(f)^m &=&\Per(\psi)\circ\Per(f)^n,\\
\Per(\varphi\circ f^m) &=&\Per(\psi\circ f^n),\\
(\varphi\circ f^m)_{|\PerPts f}&=&(\psi\circ f^n)_{|\PerPts f}.\label{eq:last}
\end{eqnarray}
By Corollary~ \ref{cor:fn-Per-f}, we may find a $k\in\NN_1$ such that
$f^k(X)\subset \PerPts f$. Then, we get from \eqref{eq:last} that
\[
\varphi\circ f^{m+k}=\psi\circ f^{n+k}
\]
which proves that $[\varphi,n]=[\psi,m]$.
This proves injectivity. To prove surjectivity take a morphism $\varphi:(X,f)\to (X',f')$ in
$\Auto(\FSet)$. Then $f,f'$ are bijections. We have
\[
  \Per'([\varphi,0])=\Per'(\Szym(\varphi))=\Per(\varphi)=\varphi_{|\PerPts f}=\varphi,
\]
which proves that $\Per$ is a surjector.
\qed

\begin{cor}
\label{cor:Szym-FSet-autom-representant}
Every object $(X,f)$ in $\Endo(FSet)$ admits an object in $\Auto(FSet)$
which is isomorphic to  $(X,f)$ in $\Szym(FSet)$. Moreover, any such object is conjugate
to $\Per(X,f)$.
\end{cor}
\proof              
It follows from Proposition~\ref{prop:f-Per} that $\Per(x,f)=(\Per f,f_{|\Per f})$ is an object in $\Auto(\FSet)$.
By Proposition~\ref{prop:universal-representant} this object is isomorphic in $\Szym(\FSet)$ to $(X,f)$.
If another object in $\Auto(\FSet)$ is isomorphic to $(X,f)$ in  $\Szym(\FSet)$ 
then it is also isomorphic to $\Per(X,f)$. Therefore, it is conjugate to $\Per(X,f)$ by Proposition~\ref{prop:Auto-conj}.
\qed


\section{Szymczak functor in $\FRel$}
\label{sec:szym_functor_in_frel}

Category $\FRel$ is the category whose objects are finite sets and whose morphisms from
set $X$ to set $X'$ consist of all relations in $X\times X'$.
The composition of morphisms $R\subset X\times X'$ and $R'\subset X'\times X''$ is defined
as the composition of relations, that is
\[
  R'\circ R:=\setof{(x,x'')\in X\times X''\mid \exists x'\in X'\; xRx' \text{ and } x'R'x''}.
\]
Then $\id_X$ is the identity morphism on $X$ for each object $X$ in $\FRel$ and
one easily verifies that so defined $\FRel$ is indeed a category.

Although the morphisms in  $\FRel$ are arbitrary relations, the following proposition shows that
isomorphisms have to be bijective maps. 

\begin{prop}
\label{prop:FRel-auto}
A relation $R\subset X\times Y$ is an isomorphism in $\FRel$ if and only if it is a bijective map.
\end{prop}
\proof
Clearly, if $R\subset X\times Y$ is a bijective map, then so is $R^{-1}$ and $R^{-1}\circ R=\id_X$ as well as  $R\circ R^{-1}=\id_Y$. Therefore, $R$  is an isomorphism in $\FRel$.
To see the converse statement assume a relation $R$ is an isomorphism. Then, there exists a relation $S\subset Y\times X$  such that $S\circ R=\id_X$ and $R\circ S=\id_Y$.
To see that $R$ is a partial map assume that $y\in R(x)$ and $y'\in R(x)$. It follows from  Proposition~\ref{prop:bijective-composition}  that $S$ is a surjective multivalued map.
Therefore, we can find a $\bar{y}\in Y$ such that $x\in S(\bar{y})$. Hence, $y\in (R\circ S)(\bar{y})=\id_Y(\bar{y})=\{\bar{y}\}$.
Similarly we get $y'\in \{\bar{y}\}$. In consequence, $y=\bar{y}=y'$ proving that $R$ is a partial map. 
It is a map, because  $X=\dom\id_X=\dom S\circ R\subset \dom R$ by Proposition~\ref{prop:dom-im-comp}.
By Proposition~\ref{prop:bijective-composition} it is a surjective map and since $X$ is finite, it is bijective map.
\qed

Given a relation  $R$  in $X$, we  set
\begin{eqnarray*}
  \gdom R&:=&\bigcap_{n\in\NN_1}\dom R^n,\\
  \gim R&:=&\bigcap_{n\in\NN_1}\im R^n,\\
  \Inv R&:=&\gdom R\,\cap\, \gim R.
\end{eqnarray*}

We say that a relation $R$ is {\em wide} if $\Inv R=X$.
We have the following proposition whose straightforward proof is left to the reader.
\begin{prop}
\label{prop:R-wide}
A relation $R$ in a finite set $X$ is wide if and only if $\dom R^n=X=\im R^n$ for all $n\in\NN_0$.
\qed
\end{prop}

Recall that a partition of a set $X$ is a family $\cA$ of mutually disjoint, nonempty subsets of $X$ such that
$X=\bigcup \cA$. Given a partition $\cA$ of $X$ and an element $x\in X$, we denote by $[x]_\cA$
the unique element of $\cA$ to which $x$ belongs.

We say that relation $R$ in $X$ is a {\em block bijection} if there exist a partition $\cA$
of $X$ and a bijection $\alpha:\cA\to\cA$ such that
\begin{equation}
\label{eq:block-bij-cA-alpha}
R=\bigcup\setof{A\times \alpha(A)\mid A\in\cA}.
\end{equation}
For a block bijection $R$ we define its size as the maximum of the cardinalities of sets in $\cA$.
Note that a bijection is always a block bijection and its size is one.

\begin{prop}
\label{prop:block-bij-struct}
Assume a relation  $R\subset X$ is a block bijection satisfying \eqref{eq:block-bij-cA-alpha}
for some partition $\cA$ of $X$ and a bijection $\alpha:\cA\to\cA$.
Then, for any $x\in X$ we have $R(x)=\alpha([x]_\cA)$.
\end{prop}
\proof
  Assume $y\in R(x)$. Then $(x,y)\in R$. It follows from \eqref{eq:block-bij-cA-alpha} that $(x,y)\in A\times \alpha(A)$
for some $A\in \cA$. Hence, $[x]_\cA=A$ and $y\in \alpha(A)=\alpha([x]_\cA)$.
To see the opposite inclusion take a $y\in \alpha([x]_\cA)$. Let $A:=\alpha([x]_\cA)$. Then $x\in A$,
$y\in \alpha(A)$ and $(x,y)\in A\times\alpha(A)\subset R$. Hence, $y\in R(x)$.
\qed

\begin{prop}
\label{prop:block-bij-bij}
Ler $R\subset X$ be a block bijection. Then, the partition $\cA$ and bijection $\alpha$ in \eqref{eq:block-bij-cA-alpha}
are uniquely determined by $R$. Moreover, if $R$ actually is a bijection, then the partition $\cA$ consists
of singletons.
\end{prop}
\proof
Let $\cA$ and $\alpha$ be the partition and bijection such that \eqref{eq:block-bij-cA-alpha} is satisfied.
It follows from Proposition~\ref{prop:block-bij-struct}
that $R(x)=\alpha([x]_\cA)$. Therefore, if  \eqref{eq:block-bij-cA-alpha} is satisfied
with $\cA$ replaced by another partition $\cB$ and $\alpha$ replaced by another bijection $\beta$,
then $\alpha([x]_\cA)=R(x)=\beta([x]_\cB)$. Since $\alpha:\cA\to\cA$ and $\beta:\cB\to\cB$ are bijections,
 this means that each set in $\cA$ equals a set in $\cB$. This is possible only if $\cA=\cB$.
 In consequence also $\alpha=\beta$. This proves the first part of the assertion.
If $R$ is a bijection, then $R(x)$ is a singleton for each $x\in X$.
It follows that $\cA$ consists of singletons.
\qed

\begin{thm}
\label{thm:part-rel}
Let $(Y,R)$ be an object in $\Endo(\FRel)$
and $(X,f)$ an object  in $\Auto(\FRel)$.
Assume that $(Y,R)$ and $(X,f)$ are isomorphic in $\Szym(\FRel)$,
that is there exist mutually inverse isomorphisms
\begin{eqnarray*}
&&[S,m]:\Szym(X,f)\to \Szym(Y,R),\\
&&[T,n]:\Szym(Y,R)\to \Szym(X,f).
\end{eqnarray*}
If $R$ is wide, then $S\circ f^k\circ T$ is a block bijection for sufficiently large $k\in\NN$
with $\setof{S(x)\mid x\in X}$ as the associated partition of $Y$.
Moreover,   $R^p$ is a block bijection for $p$ sufficiently large.
\end{thm}
\proof
Since $[S,m]$ and $[T,n]$ are mutually inverse isomorphism,
we can find a $k_0\in\NN_0$ such that 
 $T\circ S\circ f^k=f^{m+n+k}$
  and $S\circ T\circ R^k=R^{m+n+k}$ for all $k\geq k_0$.
  We will prove that
\begin{equation}
\label{eq:dom-T-im-S}
   \dom T = Y = \im S.
\end{equation}
Indeed, inclusions $\dom T \subset Y$ and $\im S \subset Y$ are obvious.
Since $R$ is wide, by Proposition~\ref{prop:R-wide} we get
$Y=\dom R^{m+n+k}$.
Hence, by Proposition~\ref{prop:dom-im-comp},  we get
\[
Y =\dom R^{m+n+k}= \dom S\circ T\circ R^k=\dom R^k\circ S\circ T\subset \dom T.
\]
Similarly, 
\[
Y=\gim R=\im  R^{m+n+k}=\im S\circ T\circ R^k\subset  \im S.
\]
This proves \eqref{eq:dom-T-im-S}.

By Proposition~\ref{prop:FRel-auto} $f$ is a bijective map. Hence, it is a wide relation and an analogous argument proves that 
\begin{equation}
\label{eq:dom-S-im-T}
   \dom S = X = \im T.
\end{equation}

Since $f$ is a bijective map, we see that $\check{f}:=f^{m+n}=T\circ S: X\to X$ is also a bijective map.
We claim that
\begin{equation}
\label{eq:S-eq-T-1}
    S(x)=T^{-1}(\check{f}(x)) \text{ for any $x\in X$.}
\end{equation}
To see this take a $y\in S(x)$. Then $x S y$. By \eqref{eq:dom-T-im-S} we may find an $x'\in X$ such that
$y T x'$. It follows that $x\,T\circ S \,x'$ which means $x'=\check{f}(x)$.
Thus, $y\in T^{-1}(x')=T^{-1}(\check{f}(x))$, which proves that $S(x)\subset T^{-1}(\check{f}(x))$.
To prove the opposite inclusion take a $y\in T^{-1}(\check{f}(x))$.
Then $y T x'$, where  $x':=\check{f}(x)$. Since $\check{f}= T\circ S$, there exists
a $y'\in Y$ such that $x S y'$ and $y' T x$. But, by \eqref{eq:dom-T-im-S}
$y\in\im S$. Therefore, we can find an $x''\in X$ such that $x'' S y$.
Hence, $x''\, T\circ S \, x'$ which means $x'=\check{f}(x'')$. It follows that
$\check{f}(x'')=\check{f}(x)$ and bijectivity of $\check{f}$ implies $x=x''$.
This together with $y\in S(x'')$ gives $y\in S(x)$ and completes the proof
of the opposite inclusion.

   We will also prove that
\begin{equation}
\label{eq:Sx-disjoint}
    S(x_1)\cap S(x_2)=\emptyset \text{ for $x_1,x_2\in X$,  $x_1\neq x_2$.}
\end{equation}
To see \eqref{eq:Sx-disjoint} assume to the contrary that
there exists a $y\in S(x_1)\cap S(x_2)$.
By \eqref{eq:dom-T-im-S} we may find an $x\in X$ such that
$x\in T(y)$. It follows that $x\in T(S(x_1))$ and $x\in T(S(x_2))$.
Since $T\circ S=\check{f}$ is a bijection, we get $x=\check{f}(x_1)$ and $x=\check{f}(x_2)$. It follows that
$x_1=\check{f}^{-1}(x)=x_2$, a contradiction proving \eqref{eq:Sx-disjoint}.

Consider the family $\cA:=\setof{S(x)\mid x\in X}$. By \eqref{eq:dom-S-im-T} the elements of $\cA$ are non-empty,
by \eqref{eq:Sx-disjoint} they are disjoint and from \eqref {eq:dom-T-im-S} we get $\bigcup\cA=Y$.
Hence, $\cA$ is a partition of $Y$. 
Fix a $k\geq k_0$ and define a bijection $\alpha:\cA\to \cA$ by $\alpha(S(x)):=S(f^k(\check{f}(x)))$.
We will prove that
\begin{equation}
\label{eq:block-bij}
S\circ f^k\circ T=\bigcup_{x\in X}S(x)\times \alpha(S(x)).
\end{equation}
Consider first a pair $(y,y')\in S\circ f^k\circ T$. Then, there exist $\bar{x},x'\in X$ such that
$y T\bar{x}$, $\bar{x}f^k x'$ and $x' S y'$. Let $x:=\check{f}^{-1}(\bar{x})$. 
It follows that $y\in T^{-1}(\bar{x})=T^{-1}(\check{f}(x))$ and, by \eqref{eq:S-eq-T-1}, $y\in S(x)$.
We also have $y'\in S(x')=S(f^k(\bar{x}))=S(f^k(\check{f}(x))=\alpha(S(x))$.
Hence $(y,y')\in S(x)\times \alpha(S(x))$, which proves that the left hand side of \eqref{eq:block-bij}
is contained in the right hand side. To prove the opposite inclusion take a pair $(y,y')\in S(x)\times \alpha(S(x))$
for some $x\in X$. Then $y\in S(x)=T^{-1}(\check{f}(x))$ which means that $y T \check{f}(x)$.
We also have $y'\in \alpha(S(x))=S(f^k(\check{f}(x)))$
or, equivalently, $\check{f}(x)\, (S\circ f^k)\, y'$.
Since $y T \check{f}(x)$, we obtain $(y,y')\in S\circ f^k\circ T$, which completes the proof of \eqref{eq:block-bij}. Therefore, $S\circ f^k\circ T$ is a block bijection.
Moreover, since $R^{m+n+k}=S\circ T\circ R^{k}=S\circ f^{k}\circ T$ holds for all sufficiently large $k$,
equation \eqref{eq:block-bij} implies that $R^p$ is a block bijection for $p$ sufficiently large.
\qed

\begin{cor}
\label{cor:Szym-rel-implies-Szym-set}
Let $(X,f)$ and $(Y,g)$ be objects in $\Endo(\FSet)$.
Then, $(X,f)$ and $(Y,g)$ are also objects in $\Endo(\FRel)$.
If objects  $(X,f)$ and $(Y,g)$
are isomorphic in $\Szym(\FRel)$, then they are also isomorphic in $\Szym(\FSet)$.
\end{cor}
\proof
  It follows from Corollary~ \ref{cor:Szym-FSet-autom-representant} that both $(X,f)$ and $(Y,g)$
are isomorphic in $\Szym(\FSet)$ to objects in $\Auto(\FSet)$.
Clearly, these isomorphisms are also isomorphisms in $\Szym(\FSet)$.
Therefore, without loss of generality we may assume that $(X,f)$ and $(Y,g)$ are objects in $\Auto(\FSet)$.
Let $[S,m]:\Szym(X,f)\to \Szym(Y,R)$ and $[T,n]:\Szym(Y,R)\to \Szym(X,f)$ be mutually inverse
isomorphisms in $\Szym(\FRel)$. Note that every bijection is obviously a wide relation.
Therefore, it follows from Theorem~ \ref{thm:part-rel}
that $R:=S\circ f^k\circ T$ is a block bijection with $\setof{S(x)\mid x\in X}$ as the associated partition of $Y$.
We also know that $S\circ T\circ g^k=g^{m+n+k}$ for a $k\in\NN$.
Hence, $S\circ f^k\circ T=S\circ T\circ g^k=g^{m+n+k}$ is a bijection.
It follows that also $S\circ T$ is a bijection.
We get from Proposition~\ref{prop:block-bij-bij} that the partition
$\setof{S(x)\mid x\in X}$  consists of singletons. This means that $S$ is a map.
It is surjective, because $\setof{S(x)\mid x\in X}$ is a partition of $Y$.
By Proposition~ \ref{prop:bijective-composition} it is also injective.
Hence, $S$ is a bijection. Since $S\circ T$ is a bijection,
it follows that also $T=S^{-1}\circ (S\circ T)$ is a bijection.
This shows that $(X,f)$ and $(Y,g)$ are conjugate.
In particular, they are isomorphic in $\Szym(\FSet)$.
\qed

\begin{prop}
\label{prop:gim-gdom-inv}
For every relation $R$ in $X$ there exists a $q\in\NN_1$ such that   for all $p\geq q$ we have 
$\gdom R=\dom R^p$ and $\gim R=\im R^p$.
\end{prop}
\proof
  Since $\dom R^n$ is a decreasing sequence of sets and $X$ is finite,
  there exists a $q\in\NN$ such that $\dom R^q = \dom R^{q+1}$.
  It follows that $\gdom R=\dom R^p$  for $p\geq q$.
  The argument for $\gim R$ is analogous.
\qed

The following proposition shows that each relation is equivalent in the Szymczak category
to a wide relation.

\begin{prop}
\label{prop:FRel-inv}
For a relation $R$ in $X$ we have
\[
\Szym(X,R)  \cong \Szym(\Inv R, R_{|\Inv R}).
\]
\end{prop}
\proof
 By Proposition~\ref{prop:gim-gdom-inv} we may fix an $n\in\NN$
 such that $\dom R^n=\gdom R$ and  $\im R^n=\gim R$.
 Let $A:=\Inv R$ and let $\bar{R}:=R_{|A}$.
 Set $S:=(R^n)_{|X\times A}$ and $T:=(R^n)_{|A\times X}$.
 We will prove that the following diagrams
\[
\begin{diagram}
 \node{X}
 \arrow{s,l}{S}
 \arrow{e,t}{R}
 \node{X}
 \arrow{s,r}{S}
 \node{X}
 \arrow{e,t}{R}
 \node{X} \\
 \node{A}
 \arrow{e,t}{\bar{R}}
 \node{A}
 \node{A}
 \arrow{n,l}{T}
 \arrow{e,t}{\bar{R}}
 \node{A}
 \arrow{n,l}{T}
\end{diagram}
\]
commute.
To see that $\bar{R}\circ S\subset S\circ R$ take $(x,y)\in \bar{R}\circ S$.
Then $x\in X$, $y\in A$ and there exists an  $a\in A$ suxh that $(x,a)\in R^n$ 
and $(a,y)\in \bar{R}\subset R$. 

Choose an $x'\in X$  such that  $(x,x')\in R$,
$(x',a)\in R^{n-1}$.
It follows that
$(x',y)\in S$ and $(x,y)\in S\circ R$.
To prove the opposite inclusion take $(x,y)\in S\circ R$.
Then, there exist $x',x''\in X$ such that  $(x,x')\in R$,
$(x',x'')\in R^{n-1}$ and $(x'',y)\in R_{|X\times A}$.
In particular, $(x,x'')\in R^n$. We will show that $x''\in A$.
Indeed, $x''\in\im R^n=\gim R$ and since $y\in A\subset\gdom R$
and $(x'',y)\in R$, it follows that $x''\in\gdom R$.
Hence,  $(x,x'')\in S$ and $(x'',y)\in\bar{R}$ which implies
$(x,y)\in \bar{R}\circ S$.
The proof of the commutativity of the other diagram is similar.

We will also prove that
\begin{eqnarray}
S\circ T&=&R^{2n}\\
T\circ S&=&\bar{R}^{2n}.
\end{eqnarray}
The inclusions $S\circ T\subset R^{2n}$ and $T\circ S\supset \bar{R}^{2n}$
follow immediately from Proposition~ \ref{prop:rel-comp-monotonicity}.
To see that $S\circ T\supset R^{2n}$ take $(x,y)\in R^{2n}$.
Then, there exists a $z\in X$ such that $(x,z)\in R^n$ and $(z,y)\in R^n$.
It follows that $z\in \im R^n=\gim R$ and $z\in \dom R^n=\gdom R$.
Hence, $z\in \Inv R=A$ and we get $(x,z)\in T$, $(z,y)\in S$ and $(x,y)\in S\circ T$.
In order to prove that $T\circ S\subset \bar{R}^{2n}$ take $(x,y)\in T\circ S$.
Then, $x,y\in A$ and there exists a sequence $x=x_0,x_1,\ldots x_n=y$ of points in $X$
such that $(x_{i-1},x_i)\in R$ for $i=1,2,\ldots 2n$. Since $x,y\in A$, it is straightforward
to observe that each $x_i\in A$. Therefore, $(x_{i-1},x_i)\in \bar{R}$, which proves
that $(x,y)\in \bar{R}^{2n}$.

Finally, we have
\begin{eqnarray}
~ [S,n]\circ[T,n]&=&[S\circ T,2n]=[R^{2n},2n]=[\id_X,0]\\
~ [T,n]\circ[S,n]&=&[T\circ S,2n]=[\bar{R}^{2n},2n]=[\id_A,0],
\end{eqnarray}
which proves that $[S,n]:\Szym(X,R)\to\Szym(A,\bar{R})$ and $[T,n]:\Szym(A,\bar{R})\to\Szym(X,R)$
are mutually inverse isomorphisms.
\qed

\begin{prop}
\label{prop:eventual-period}
Let $(X,R)$ be an object of $\Endo(\FRel)$.
Then there exists a $p\in\NN_1$ such that
\begin{equation}
\label{eq:eventual-period-1}
     R^{i+p}=R^i \text{ for } i\geq p
\end{equation}
and, in particular,
\begin{equation}
\label{eq:eventual-period-2}
     R^{kp}=R^p \text{ for } k\in\NN_1.
\end{equation}
\end{prop}


\proof
Since $X$ is finite, the set of all relations in $X$ is finite.
In particular, the set of values of the sequence  $\NN_1\ni R\mapsto R^n$ is finite.
It follows that there exist $m_1,m_2\in\NN_1$ such that $m_1<m_2$ and $R^{m_1}=R^{m_2}$.
Set $q:=m_2-m_1$ and choose an $m\in\NN_1$ such that $p:=mq\geq m_1$.
Then $R^{m_1+q}=R^{m_1}$. Multiplying both sides by $R^{q}$
we obtain $R^{m_1+2q}=R^{m_1+q}=R^{m_1}$. Thus, an induction argument proves that $R^{m_1+kq}=R^{m_1}$
for $k\in\NN_1$.
Fix $i\geq p$. Then $i\geq m_1$ and
\[
   R^{i+p}=R^{(i-m_1)+m_1+mq}=R^{(i-m_1)+m_1}=R^i,
\]
which proves \eqref{eq:eventual-period-1} and \eqref{eq:eventual-period-2}
follows easily from \eqref{eq:eventual-period-1} by induction.
\qed

We refer to a $p$ satisfying \eqref{eq:eventual-period-1} as an {\em eventual period} of $R$. The key feature of an eventual period is \eqref{eq:eventual-period-1}. Therefore, we do not require that the eventual period be the smallest number with this property.
Note that a similar concept, named index, is introduced in \cite{KR1986}.

\begin{thm}
\label{thm:eventual-period}
Let $(X,R)$ be an object of $\Endo(\FRel)$ and let $p$ be an eventual period of $R$.
Then for each $s\in\NN_1$
\[
\Szym(X,R^s)\cong\Szym(X,R^{s+p}).
\]
\end{thm}
\proof
   Let $S:=T:=R^p$. We claim that $[S,p]:\Szym(X,R^s)\to\Szym(X,R^{s+p})$
   and $[T,p]:\Szym(X,R^{s+p})\to \Szym(X,R^s)$ are mutually inverse isomorphisms in $\Szym(\FRel)$.
   Since $p+s\geq p$, we get from \eqref{eq:eventual-period-1} that
\[
   R^{p+s}\circ T=  R^{2p+s} =  R^{p+s}  = T\circ R^s,
\]
\[
   R^{s}\circ S=  R^{p+s} =  R^{2p+s}  = S\circ R^{p+s}.
\]
This shows that $R$ and $S$ are morphisms in $\Endo(\FRel)$. Moreover,
by \eqref{eq:eventual-period-2}
\[
   T\circ S=  R^{2p} = R^p= R^{2sp}=(R^s)^{2p},
\]
\[
   S\circ T=  R^{2p}  = R^p= R^{2(s+p)p}=(R^{s+p})^{2p},
\]
which proves that $[T,p]\circ [S,p]=[\id_X,0]$ and  $[S,p]\circ [T,p]=[\id_X,0]$, that is
$[T,p]$ and $[S,p]$ are mutually inverse isomorphisms.
\qed

The following proposition is straightforward. 
\begin{prop}
\label{prop:connected-is-wide}
Let $R\subset X\times X$ be a strongly connected relation. Then $R$ is wide. Moreover, if $x\in R^k(x)$, then $x\in R^{kl}(x)$ for each $l\in\NN_1$. \qed
\end{prop}

Let $\gcd(a,b)$ denote the greatest common divisor of $a,b\in\ZZ$. Consider the greatest common divisor of the length of all cycles in a strongly connected relation $R$. We call this number the {\em period of $R$}. In order to compute the period of $R$ one can consider the set of cycles with different vertices (cf. \cite[Definition 7.1.]{St1994}). The following proposition relates the period of a strongly connected relation with its eventual period.
\begin{prop}
\label{prop:period-and-ev-period}
Let $p\in\NN_1$ be an eventual period of a strongly connected relation $R\subset X\times X$ and let  $q\in\NN_1$ be the period of $R$. Then $q\leq p$. Moreover, $q|p$.
\end{prop}
\proof
Assume to the contrary that $q>p$. Then there exists at least one $x\in X$ such that $x\notin R^p(x)$ because otherwise $q$ would divide $p$. Since $R$ is strongly connected, there exists an $l\in\NN_1$ such that $x\in R^l(x)$ and $q|l$. By Proposition \ref{prop:eventual-period} we get $R^{lp}(x)=R^p(x)\not\ni x$. It follows from Proposition \ref{prop:connected-is-wide} that $x\in R^{lp}(x)$, a contradiction.

In order to prove that $q|p$ note that for any $x\in X$ there exists an $i\in\{0,\ldots,p-1\}$ such that $x\in R^{p+i}(x)$. Indeed, by Proposition \ref{prop:eventual-period} for any $m\geq p$ we have $R^m=R^{p+i}$ for some $i\in\{0,1,\ldots,p-1\}$. From the same proposition we conclude that for any $k\in\NN_1$ the equation $R^{p+i}=R^{kp+i}$ holds. Therefore, $x\in R^{kp+i}(x)$ and this means $q|p+i$ and $q|kp+i$. It follows that $q|a(p+i)+b(kp+i)$ for any $a,b\in\ZZ$ and $k\in\NN_1$. Setting $a=-1$, $b=1$ and $k=2$ we get $q|p$.
\qed

There are some relationships between eventual periods of an arbitrary relation and eventual periods of the relation restricted to its strongly connected components. 

\begin{prop}
\label{prop:powers-of-restrictions}
Let $U\subset X$ be a strongly connected component of an arbitrary $R\subset X\times X$. Then
$$
(R|_U)^n=(R^n)|_U
$$
for each $n\in\NN$.
\end{prop}
\proof
The left-hand-side is clearly contained in the right-hand-side.

To prove the opposite inclusion consider a pair $(u,v)\in(R^n)|_U$. Then $(u,v)\in U\times U$ and there is a $(u,v)$-walk in $R$ of length $n$. Since $u$ and $v$ belong to the same strongly connected component of $R$, there is a $(v,u)$-walk in ${R|_U}$. Concatenation of both walks gives a cycle in ${R|_U}$, because $U$ is a strongly connected component of $R$. Therefore, vertices lying on the $(u,v)$-walk belong to $U$. In consequence, $(u,v)\in (R|_U)^n$.
\qed

\begin{cor}
\label{cor:comp-period-divides-ev-period}
Let $U\subset X$ be a strongly connected component of $R\subset X\times X$ and let $p$ and $p_U$ be eventual periods of $R$ and $R|_U$, respectively. Then $p_U|p$.
\end{cor}
\proof
We have $R^{p+i} = R^i$ for $i\geq p$. By Proposition \ref{prop:powers-of-restrictions}
$$
(R|_U)^{p+i} = (R^{p+i})|_U=(R^i)|_U=(R|_U)^i.
$$
Hence, $p$ is a multiple of $p_U$.
\qed

\begin{prop}
\label{prop:walk-length-divisible-by-period}
Let $q\in\NN_1$ be the period of a strongly connected relation $R\subset X\times X$ and let $x,y\in X$. Then the following conditions are equivalent:
\begin{itemize}
\item[(i)] there exists an $(x,y)$-walk in $R$ with length divisible by $q$,
\item[(ii)] each $(x,y)$-walk in $R$ has length divisible by $q$.
\end{itemize}
\end{prop}
\proof
Let $c=x\ldots y$ be an $(x,y)$-walk in $R$ such that $q|\#c$. Consider a walk $d=x\ldots y$ in $R$ such that $\#c\neq\#d$. Since $R$ is strongly connected, there exists a $(y,x)$-walk $e$ in $R$. Then $ce$ is a cycle passing through the vertex $y$. Since $q$ is the period of $R$, $q$ divides the length of the cycle. Also $q|\#c$, hence $q|\#e$. Since $de$ is also a cycle in $R$ passing through $y$, the period $q$ divides its length. Therefore, $q|\#d$, because $q|\#e$.

To prove the opposite implication it suffices to note that the existence of an $(x,y)$-walk follows from the strong connectivity of $R$.
\qed

Let $R\subset X\times X$ be an arbitrary relation. We write $x\rightarrow_R y$ to denote that there is a walk in $R$ from $x$ to $y$ of positive length. We say that $x,y\in X$ are {\em strongly connected} and write $x\leftrightarrow_R y$ if $x\rightarrow_R y$ and $y\rightarrow_R x$. Relation $\leftrightarrow_R$ is clearly symmetric and transitive. Hence, it is an equivalence relation in 
$$
X_R:=\{x\in X\ |\ x\leftrightarrow_R x\}.
$$
We call $X_R$ the {\em recurrent set} of $R$ and its elements the {\em recurrent vertices} of $R$. The equivalence classes of $\leftrightarrow_R$ in $X_R$ are called the {\em strongly connected components} of $R$. For a recurrent vertex $x\in X_R$ we denote by $[x]_R$ the strongly connected component to which $x$ belongs.

We refine the relation $\leftrightarrow_R$ in $X_R$ to a relation $\sim_R$ in $X_R$ by defining $x\sim_R y$ for $x,y\in X_R$ if $x\leftrightarrow_R y$ and each walk from $x$ to $y$ has length equal to zero modulo the period of $R$ restricted to the common strongly connected component of $x$ and $y$.


Notice that if $R\subset X\times X$ is a strongly connected relation, then $X_R = X$ and $\leftrightarrow_R$ has exactly one equivalence class.

\begin{prop}[{cf. \cite[Lemma 6, Corollary 1]{JS1996}}]
\label{prop:period-rel-is-equivalence-rel}
Let $q\in\NN_1$ be the period of a strongly connected relation $R\subset X\times X$. Then $\sim_R$ is an equivalence relation in $X$ with exactly $q$ distinct equivalence classes. \qed
\end{prop}

\begin{cor}
\label{cor:sim-is-equivalence-rel-in-recurrent-set}
Let $R\subset X\times X$ be an arbitrary relation. Then $\sim_R$ is an equivalence relation in $X_R$. \qed
\end{cor}

In order to proceed we need the following fact about the existence of solutions of particular Diophantine equation. \footnote{Proposition \ref{prop:positive-solutions} and the proof after D. Jao from \\ {https://djao.math.uwaterloo.ca/w/Positive\_Solutions\_to\_linear\_diophantine\_equation.} (Retrieved February 4, 2021)}

\begin{prop}[D. Jao]
\label{prop:positive-solutions}
Let $a,b\in\NN_1$ and let $q:=\gcd(a,b)$. For any $n\in\NN$ the equation
$$
ax+by=\frac{ab}{q}+q+nq
$$
has a solution $x,y$ consisting of positive integers.
\end{prop}
\proof
Let $a':=\frac{a}{q}\in\NN$ and let $b':=\frac{b}{q}\in\NN$. Then $a'$ and $b'$ are coprime. We will show that 
\begin{equation}
\label{eq:reminders-set}
\{b'y \mkern-12mu\mod a'\ |\ y=1,2,\ldots,a'\}=\{0,1,\ldots,a'-1\}.
\end{equation}
Clearly, the left-hand-side of (\ref{eq:reminders-set}) is contained in the right-hand-side. Assume there is no equality.  Then there exist $y_1,y_2\in\{1,\ldots,a'\}$ such that $b'y_1 = b'y_2 \mkern-8mu\mod a'$ and $y_1<y_2$. It follows that $a'$ divides $b'y_2-b'y_1$ and there exists an $s\in\NN$ such that $a's=b'(y_2-y_1)$. Since $a',b'$ are coprime, $b'$ divides $s$ and $u:=\frac{s}{b'}\in\NN$. We get $a'\leq a'u = y_2-y_1<a'$, a contradiction which proves (\ref{eq:reminders-set}).

Let $r$ denote the reminder of the division of $1+n$ by $a'$. Then $1+n=ka'+r$. By (\ref{eq:reminders-set}) we may choose a $y\in\{1,\ldots,a'\}$ such that $b'y \mkern-8mu\mod a' = r$. Then $a'$ divides $1+n-b'y$ and, in consequence, $a$ divides $q(1+n-b'y)+b'a = \frac{ab}{q}+q+nq-by$. Hence,
$$
x:=\frac{\frac{ab}{q}+q+nq-by}{a}
$$
is an integer. It remains to show that $x>0$. For this end note that $by\leq ba'$ and
$$
ax\geq ba'+q+nq-ba' = q+nq.
$$
Therefore $x\geq \frac{q}{a} + \frac{nq}{a} > 0$.
\qed

\begin{lem}
\label{lem:connecting-cycle}
Let $R\subset X\times X$ be a strongly connected relation with its period equal to $q$ and an eventual period equal to $p$. Then for all cycles $c$ and $d$ in $R$ of length equal to $k_1$ and $k_2$ respectively, there exists a cycle $b$ of length $\gcd(k_1,k_2)+p\gcd(k_1,k_2)$ such that each vertex of cycles $c$ and $d$ lies on cycle $b$.
\end{lem}
\proof
Assume that a vertex $x$ lies on a cycle $c$ (for a vertex lying on $d$ the proof is analogous). Then $x\in R^{k_1}(x)$. In general, cycles $c$ and $d$ need not share vertices. By connectivity of $R$, there exist a walk from $x$ to any vertex $y$ lying on $d$ as well as a walk from $y$ to $x$. The concatenation of the $(x,y)$-walk with the $(y,x)$-walk is a cycle $e$ in $R$ with $x$ and $y$ lying on $e$. Let us assume that $\#e:=l$.

By Proposition \ref{prop:positive-solutions}, for each $m\in\NN$ there exists a solution $l_1,l_2\in\NN_1$ of the equation
$$
k_1l_1+k_2l_2 = \frac{k_1k_2}{\gcd(k_1,k_2)} + \gcd(k_1,k_2) + m\gcd(k_1,k_2).
$$
We have $k_1':=\frac{k_1}{\gcd(k_1,k_2)}\in\NN$. Hence,
$$
k_1l_1+k_2l_2 = \gcd(k_1,k_2)(\frac{k_1'k_2}{\gcd(k_1,k_2)} + 1 + m).
$$
Putting $m:= m_1p-\frac{k_1'k_2}{\gcd(k_1,k_2)}$ where $m_1\in\NN$ is large enough to ensure that $m>0$, we get a solution $l_1,l_2\in\NN_1$ such that
$$
k_1l_1+k_2l_2 = \gcd(k_1,k_2) + m_1p\gcd(k_1,k_2).
$$
It follows that $x\in R^{k_1l_1+k_2l_2+lp\gcd(k_1,k_2)}(x)$, because, by Proposition \ref{prop:connected-is-wide}, cycle $e$ may be concatenated with itself $p\gcd(k_1,k_2)$ times, and cycles $c$ and $d$ may be concatenated with themselves $l_1$ and $l_2$ times respectively. Hence, by Proposition \ref{prop:eventual-period},
$$
x\in R^{k_1l_1+k_2l_2+lp\gcd(k_1,k_2)}(x) = R^{\gcd(k_1,k_2)+p\gcd(k_1,k_2)}(x).
$$
This means that there exists a cycle in $R$ with $x$ lying on that cycle of length $\gcd(k_1,k_2)+p\gcd(k_1,k_2)$. 
\qed

\begin{prop}
\label{prop:identity-in-p-plus-q-power}
Let $R\subset X\times X$ be a strongly connected relation with its period equal to $q$. For every eventual period $p$ of $R$ we have $\id_X\subset R^{p+q}$.
\end{prop}
\proof
Let $x\in X$ and let $\{k_1,\ldots,k_n\}$ be the set of lengths of all cycles with different vertices in $R$. We have $q=\gcd(k_1,\ldots,k_n)$ (see \cite[Definition 7.1]{St1994}). Take a cycle such that $x$ lies on it. Assume that the cycle has length $k_1$. Hence, $x\in R^{k_1}(x)$. Take the next cycle of length  $k_2$. By Lemma \ref{lem:connecting-cycle}, there exists a cycle of length $\gcd(k_1,k_2)+p\gcd(k_1,k_2)$ with $x$ lying on it. 

Let us take the next cycle of length $k_3$. Note that
$$
x\in R^{\gcd(k_1,k_2)+p\gcd(k_1,k_2)}(x) = R^{\gcd(k_1,k_2)+k_3p\gcd(k_1,k_2)}(x),
$$
and by the identity $\gcd(a+cb,c)=\gcd(a,c)$ for $a,b,c\in\ZZ$ we get
$$
\gcd(\gcd(k_1,k_2)+k_3p\gcd(k_1,k_2),k_3) = \gcd(\gcd(k_1,k_2),k_3) = \gcd(k_1,k_2,k_3).
$$
Applying again Lemma \ref{lem:connecting-cycle} for cycles of lengths $\gcd(k_1,k_2)+p\gcd(k_1,k_2)$ and $k_3$, we can find a cycle such that
$$
x\in R^{\gcd(k_1,k_2,k_3) + p\gcd(k_1,k_2,k_3)}(x).
$$ 
Continuing for the remaining $k_i$ we get $x\in R^{q+pq}(x) = R^{q+p}(x)$ which implies $\id_X\subset R^{p+q}$.
\qed

As a corollary to the above proposition we get a variant of Proposition \ref{prop:eventual-period} for strongly connected relations.

\begin{cor}
\label{cor:strong-relation-after-ev-period}
Let $R$ be a strongly connected relation with its period equal to $q$ and an eventual period equal to $p$. Then
\begin{equation}
\label{eq:powers-after-ev-period}
R^{p+kq} = R^{p} \text{ for } k\in\NN.
\end{equation}
\end{cor}
\proof
We prove inductively on $k\in\NN$ that
\begin{equation}
\label{eq:rel-powers-inclusion}
R^{p+kq}\subset R^{p+(k+1)q}.
\end{equation}
By Proposition \ref{prop:identity-in-p-plus-q-power} we have $\id_X\subset R^{p+q}$, hence $R^p\subset R^{2p+q}$. By Proposition \ref{prop:eventual-period} we have $R^{2p+q}=R^{p+q}$. This proves (\ref{eq:rel-powers-inclusion}) for $k=0$. 

Proceeding by induction we get
$$
R^{p+(k+1)q} = R^{p+kq}\circ R^{q} \subset R^{p+(k+1)q}\circ R^q = R^{p+(k+2)q},
$$
which completes the proof of (\ref{eq:rel-powers-inclusion}).

We will now prove (\ref{eq:powers-after-ev-period}). By Proposition \ref{prop:period-and-ev-period},  $p=mq$ holds for some $m\in\NN_1$. Fix an $s\in\NN$ such that $sm\geq k$. By (\ref{eq:rel-powers-inclusion}), we have
\[
\pushQED{\qed}
R^p\subset R^{p+kq}\subset R^{p+smq} = R^{p+sp} = R^p. \qedhere
\popQED
\]

We are now ready to present a theorem expressing the equivalence classes of $\sim_R$ in $X_R$ in terms of a power of relation $R\subset X\times X$.


\begin{thm}
\label{thm:partition}
Let $R\subset X\times X$ be an arbitrary relation and let $p\in\NN_1$ be an eventual period of $R$. Then for each $x\in X_R$ we have
\begin{equation}
\label{eq:period-classes}
[x]_{\sim_R} = R^p(x)\cap [x]_R.
\end{equation}
In particular, if $R$ is a strongly connected relation, then $[x]_{\sim_R}=R^p(x)$.
\end{thm}
\proof
Let $y\in [x]_{\sim_R}$. This means that there exists an $(x,y)$-walk of length $kq$, where $q\in\NN_1$ is the period of $R|_{[x]_R}$ and $k\in\NN_1$. In other words, $y\in (R|_{[x]_R})^{kq}(x)$. Notice that $x\in (R|_{[x]_R})^p$. Indeed, we have 
$$
x\in\id_{[x]_R}(x)\subset (R|_{[x]_R})^{p_{[x]_R}+q}(x)=(R|_{[x]_R})^{p_{[x]_R}}(x)\subset (R|_{[x]_R})^{p}(x),
$$
where $p_{[x]_R}$ is an eventual period of $R|_{[x]_R}$. 
By Proposition \ref{prop:identity-in-p-plus-q-power}, Corollary \ref{cor:strong-relation-after-ev-period} and Proposition \ref{prop:powers-of-restrictions} we get 
$$
y\in (R|_{[x]_R})^{p+kq}(x) = (R|_{[x]_R})^{p}(x) = (R^p)|_{[x]_R}(x)\subset R^p(x).
$$ 
It is clear that $y\in [x]_R$.

In order to prove the opposite inclusion take a $y\in R^p(x)\cap [x]_R$. There exists an $(x,y)$-walk of length $p$ in $R|_{[x]_R}$. Since $R|_{[x]_R}$ is strongly connected, there exists also a $(y,x)$-walk of length $l$ in $R|_{[x]_R}$ for some $l\in\NN_1$. Concatenation of these walks is a cycle of length $p+l$. Hence, the period $q$ of  $R|_{[x]_R}$ divides $p+l$. By Proposition \ref{prop:period-and-ev-period}, we have $q|p_{[x]_R}$, where $p_{[x]_R}$ is an eventual period of $R|_{[x]_R}$. By Corollary \ref{cor:comp-period-divides-ev-period}, $q|p$. Therefore, $q|l$ and this proves $y\sim_R x$, that is $y\in [x]_{\sim_R}$.
\qed

Let $(X,R)\in\Endo(\FRel)$ and let $p\in\NN_1$ be an eventual period of $R$. Since $\sim_R$ is an equivalence relation in $X_R$, we may consider the object $(X_R/_{\sim_R},\bar{R})\in\Endo(\FRel)$, where $\bar{R}$ is induced on equivalence classes of $\sim_R$ given by
\begin{equation}
\label{eq:R-bar}
[x]_{\sim_R} \bar{R}\ [y]_{\sim_R} \ \text{ if }\ x R^{p+1} y
\end{equation}
for $x,y\in X_R$. The relation $\bar{R}$ is well-defined. This is a consequence of the following implication:
$$
x\sim_R x',\ x R^{p+1} y,\ y\sim_R y' \implies x' R^{p+1} y'.
$$
The implication holds. Indeed, there are an $(x',x)$-walk and a $(y,y')$-walk of length equal to zero modulo the period of the strongly connected component containing $x,x'$ and $y,y'$, respectively. By Corollary \ref{cor:strong-relation-after-ev-period} and Corollary \ref{cor:comp-period-divides-ev-period} there are also an $(x',x)$-walk and a $(y,y')$-walk of length equal to $p$. Concatenating these walks of length equal to $p$ with an $(x,y)$-walk of length equal to $p+1$ in the right order we get the $(x',y')$-walk of length equal to $3p+1$. By Proposition \ref{prop:eventual-period}, there is an $(x',y')$-walk of length equal to $p+1$ which proves the implication.

\begin{lem}
\label{lem:R-bar-property}
Let $(X,R)\in\Endo(\FRel)$ and let $p\in\NN_1$ be an eventual period of $R$. Then for $\bar{R}$ given by (\ref{eq:R-bar}) 
$$
\bar{R}([x]_{\sim_R}) = \bar{R}(\{[y]_{\sim_R}\ |\ y\in R^p(x), \ y\in X_R\})
$$
for all $x\in X_R$. Moreover, $p$ is an eventual period of $\bar{R}$.
\end{lem}
\proof
The left-hand-side is clearly contained in the right-hand-side.

To prove the opposite inclusion consider a $[z]_{\sim_R}$ which belongs to the right-hand-side. This means that there is a $y\in R^p(x)$, $y\in X_R$ such that $[y]_{\sim_R}\bar{R} [z]_{\sim_R}$. Thus, $y R^{p+1} z$ and $z\in R^{p+1}(y)\subset R^{p+1}(R^p(x)) = R^{p+1}(x)$. It follows that $x R^{p+1} z$ and $[x]_{\sim_R}\bar{R} [z]_{\sim_R}$.

Let $i\geq p$. We have
$$
\begin{array}{ccl}
\bar{R}^i([x]_{\sim_R}) & = & \{[y]_{\sim_R}\ |\ x (R^{p+1})^i y\} = \{[y]_{\sim_R}\ |\ x R^{pi+i+p^2+p} y\} \\
 & = &\{[y]_{\sim_R}\ |\ x (R^{p+1})^{i+p} y\} = \bar{R}^{i+p}([x]_{\sim_R}), \\
\end{array}
$$
which proves that $p$ is an eventual period of $\bar{R}$.
\qed

\begin{lem}
\label{lem:sccv}
Let $R\subset X\times X$ be an arbitrary relation and let $p\in\NN_1$ be an eventual period of $R$. For each $x\in X$ and $n\in\NN$
\begin{equation}
\label{eq:relation-on-components}
R^{p+n}(x) = R^p(R^{p+n}(x)\cap X_R).
\end{equation}
\end{lem}
\proof
Note that if $X_R=\emptyset$, then the relation $R^p$ is empty. Therefore, in this case the theorem is trivial. Hence, assume that $X_R\neq\emptyset$. We prove formula (\ref{eq:relation-on-components}) inductively on $n\in\NN$. 

Assume that $n=0$. We need to prove that $R^{p}(x) = R^p(R^{p}(x)\cap X_R)$ for each $x\in X$. For the proof of the right-to-left inclusion, note that for each $x\in X$ we have $R^p(x)\cap X_R\subset R^p(x)$ and, in consequence, $R^p(R^p(x)\cap X_R)\subset R^{p+p}(x)=R^p(x)$.

In order to prove the opposite inclusion take a $y\in R^p(x)$. We claim that there is an $(x,y)$-walk in $R$ such that there exists a $z\in X_R$ which belongs to the walk. Indeed, if this were not true, then we would get a contradiction in the equality (\ref{eq:eventual-period-2}) in Proposition \ref{prop:eventual-period}, because from the finiteness of $X$ there would be a number $k\in\NN_1$ such that $R^i(x)=\emptyset$ for each $i\geq k$, in particular $y\in R^p(x)=R^{kp}(x)=\emptyset$.

Let us take a $z\in U$ lying on the $(x,y)$-walk in some strongly connected component $U$. Assume that $z\in R^l(x)$ for some $l\in\{0,\ldots,p\}$. Clearly, $y\in R^{p-l}(z)$. Note that there exists a $z'\in U\subset X_R$ lying on a cycle starting at $z$ of length $p$ such that $z'\in R^{p-l}(z)$. Note that it may happen that $z=z'$. By Theorem \ref{thm:partition}, the set $R^p(z')$ contains some equivalence class of $\sim_R$ defined in $X_R$. In particular,  $z'\in R^p(z')$. Therefore,
$$
z'\in R^p(z')\subset R^{2p-l}(z)\subset R^{2p}(x) = R^p(x).
$$
In consequence, $z'\in R^p(x)\cap X_R$. We will show that there is a $(z',z)$-walk in $R$ of length $p+l$. Indeed, from the definition of $z'$ we know that $z\in R^l(z')$. Together with $z'\in R^p(z')$, we get $z\in R^{p+l}(z')$. We have $y\in R^{p-l}(z)\subset R^{2p}(z') = R^p(z')$ which proves $y\in R^p(z')\subset R^p(R^p(x)\cap X_R)$.

Hence, formula (\ref{eq:relation-on-components}) for $n=0$ is proved. Now assume that (\ref{eq:relation-on-components}) holds. We prove that (\ref{eq:relation-on-components}) also holds with $n$ replaced by $n+1$. Using the inductive assumption and the formula that the image of a union under a multivalued map is equal to the union of the images, we get
$$
\begin{array}{rcl}
R^p(R^{p+n+1}(x)\cap X_R) & = & R^p(R^{p+n}(R(x))\cap X_R) \\
 & = & R^p((\bigcup_{t\in R(x)}R^{p+n}(t))\cap X_R) \\
 & = & R^p(\bigcup_{t\in R(x)}R^{p+n}(t)\cap X_R) \\
 & = & \bigcup_{t\in R(x)}R^p(R^{p+n}(t)\cap X_R) \\
 & = & \bigcup_{t\in R(x)} R^{p+n}(t) \\
 & = & R^{p+n+1}(x),\\
\end{array} 
$$
which ends the proof.
\qed

\begin{lem}
\label{lem:eq-rel-property}
Let $R\subset X\times X$ be an arbitrary relation and let $p\in\NN_1$ be an eventual period of $R$. Then 
$$
x\sim_R x' \implies R^p(x)=R^p(x').
$$
\end{lem}
\proof
Let $x\sim_R x'$. By Theorem \ref{thm:partition}, we have $x'\in R^p(x)$ and, in consequence, $R^p(x')\subset R^p(x)$. 

The right-to-left inclusion follows by symmetry of $\sim_R$.
\qed

\begin{thm}
\label{thm:arisom}
Let $(X,R)\in \Endo(\FRel)$ and let $p\in\NN_1$ be an eventual period of $R$. Then 
$$
\Szym(X,R)\cong \Szym (X_R/_{\sim_R},\bar{R}),
$$
where $\bar{R}$ is induced on equivalence classes of $\sim_R$ given by (\ref{eq:R-bar}).
\end{thm}
\proof
Consider relations $S\subset X\times X_R/_{\sim_R}$ and $T\subset X_R/_{\sim_R}\times X$ defined by $S(x):=\{[y]_{\sim_R}\ |\ y\in R^p(x), \ y\in X_R\}$ for $x\in X$ and $T([x]_{\sim_R}):=R^p(x)$ for $[x]_{\sim_R}\in X_R/_{\sim_R}$. By Lemma \ref{lem:eq-rel-property}, $T$ is well-defined. We claim that $S$ and $T$ are morphisms in $\Endo(\FRel)$. Note that by Lemma \ref{lem:R-bar-property}, for $x\in X$ we have
$$
\begin{array}{ccl}
(S\circ R)(x) &=& S(R(x)) = \{[y]_{\sim_R}\ |\ y\in R^{p+1}(x), \ y\in X_R\} \\
 &=& \{[y]_{\sim_R}\ |\ [x]_{\sim_R}\bar{R} [y]_{\sim_R}\} = \bar{R}([x]_{\sim_R}) \\
 &=& \bar{R}(\{[y]_{\sim_R}\ |\ y\in R^p(x), \ y\in X_R\}) = (\bar{R}\circ S)(x), \\
\end{array}
$$
and, by Lemma \ref{lem:sccv}, for $[x]_{\sim_R}\in X_R/_{\sim_R}$
$$
\begin{array}{ccl}
(R\circ T)([x]_{\sim_R}) &=& R(R^p(x)) = R^{2p+1}(x) \\
&=& R^p(\{y\ |\ y\in R^{p+1}(x), y\in X_R\}) \\
&=& T(\{[y]_{\sim_R}\ |\ [x]_{\sim_R}\bar{R} [y]_{\sim_R}\}) = (T\circ \bar{R})([x]_{\sim_R}), \\
\end{array}
$$
which proves that $S$ and $T$ are morphisms in $\Endo(\FRel)$. 

Now we prove that  
$$
[S,p]\colon\Szym(X,R)\to\Szym(X_R/_{\sim_R},\bar{R})
$$ and 
$$
[T,p]\colon\Szym(X_R/_{\sim_R},\bar{R})\to\Szym(X,R)
$$ 
are mutually inverse isomorphisms in $\Szym(\FRel)$. Again by Lemma \ref{lem:sccv}, for $x\in X$ we get 
$$
(T\circ S)(x) = T(\{[y]_{\sim_R}\ |\ y\in R^p(x)\cap X_R\}) = R^p(\{y\ |\ y\in R^p(x)\cap X_R\}) = R^p(x)
$$
and for $[x]_{\sim_R}\in X_R/_{\sim_R}$
$$
\begin{array}{ccl}
(S\circ T)([x]_{\sim_R}) &=& S(R^p(x)) = \{[y]_{\sim_R}\ |\ y\in R^p(R^p(x)), y\in X_R\} \\
 &=& \{[y]_{\sim_R}\ |\ y\in R^p(x)\cap X_R\} \\
 &=& \{[y]_{\sim_R}\ |\ y\in(R^{p+1})^p(x)\cap X_R\} \\
  &=& \{[y]_{\sim_R}\ |\ [x]_{\sim_R}\bar{R}^p [y]_{\sim_R}\} = \bar{R}^p([x]_{\sim_R}). \\
\end{array}
$$
Note that, in particular, the following holds:
$$
\id_X\circ R^{2p+p} = R^{p+p} = R^p\circ R^p.
$$
Hence, $[T,p]\circ[S,p] = [T\circ S,2p] = [R^p,2p] = [\id_X,0]$. By Lemma \ref{lem:R-bar-property}, we get
$$
\id_{X_R/_{\sim_R}}\circ\bar{R}^{2p+p} = \bar{R}^p\circ \bar{R}^p. 
$$
Therefore, $[S,p]\circ[T,p] = [S\circ T,2p] = [\bar{R}^p,2p] = [\id_{X_R/_{\sim_R}},0]$, which ends the proof.
\qed

\begin{figure}[h]
\begin{center}
  \includegraphics[width=0.85\textwidth]{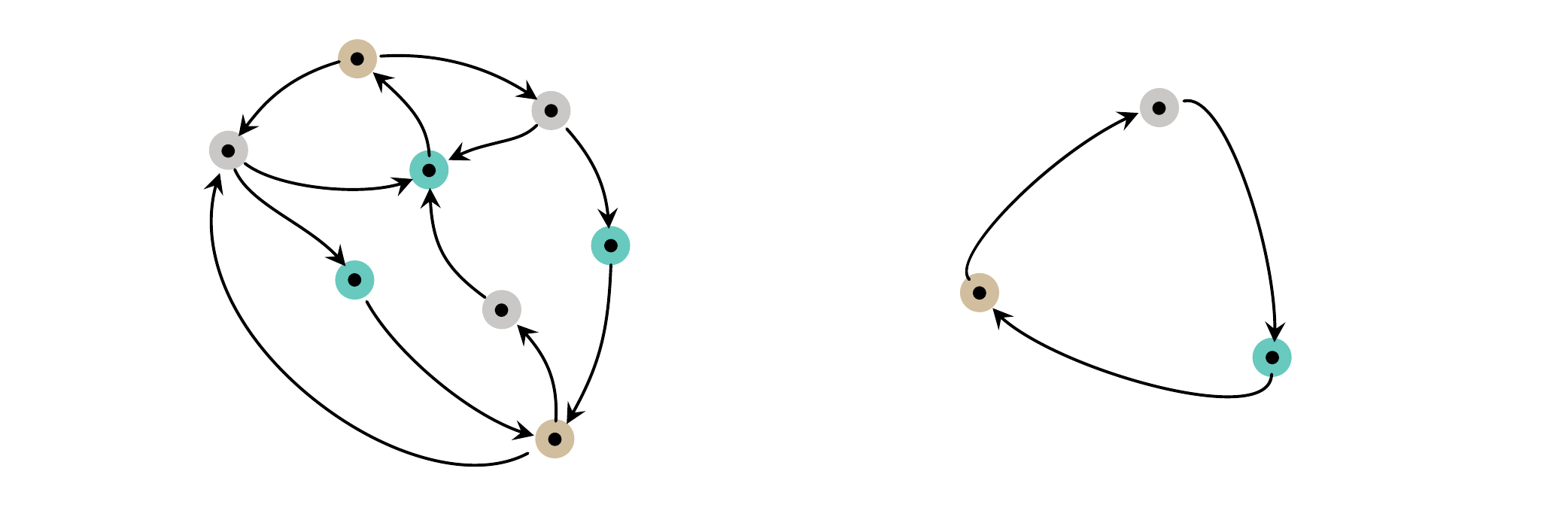}
\end{center}
  \caption{Two isomorphic relations in Szymczak's category. The eventual period $p$ and the period $q$ of the relation on the left are both $p=q=3$. The equivalence classes of the relation $\sim_R$ are marked with colors.}
  \label{fig:mv-rel}
\end{figure}

Note that for a strongly connected relation $R$, the relation $\bar{R}$ from Theorem \ref{thm:arisom} is, in fact, a cyclic bijection (see the example in Figure \ref{fig:mv-rel}).

The example in Figure \ref{fig:mv-rel} shows the partition of the set of vertices into the equivalence classes of the relation $\sim_R$.

%
\subsection{Objects in canonical form}
\label{subsec:obj-in-canonical-form}

Now we will consider particular class of objects in $\Endo(\FRel)$. We say that $(X,R)\in\Endo(\FRel)$ is in {\em canonical form} if the following conditions apply:
\begin{itemize}
\item[(i)] $X=X_R$; in other words, each element of $X$ belongs to a strongly connected component of $R$,
\item[(ii)] for each $x\in X$ the restriction $R|_{[x]_R}$ is a bijection, 
\item[(iii)] for each $n\in\NN_1$ the equation $R^{n+p}=R^n$ holds, where $p\in\NN_1$ is an eventual period of $R$.  
\end{itemize}
Note that the condition (iii) implies that the bijection from (ii) is cyclic. 

\begin{thm}[Theorem \ref{thm:main1}]
\label{thm:each-obj-has-its-canonical-form}
For each $(X,R)\in\Endo(\FRel)$ there exists an object $(\bar{X},\bar{R})\in\Endo(\FRel)$ in canonical form such that 
$$
\Szym(X,R)\cong\Szym(\bar{X},\bar{R}).
$$
\end{thm}
\proof
Let $p\in\NN_1$ be an eventual period of $R$. Consider $(\bar{X},\bar{R})$, where $\bar{X}:=X_R/_{\sim_R}$ and $\bar{R}$ is induced by $R$ on equivalence classes of $\sim_R$ as in (\ref{eq:R-bar}). We claim that $(\bar{X},\bar{R})$ is in canonical form. 

To prove that $\bar{X} = \bar{X}_{\bar{R}}$ let $\alpha\in\bar{X}$ and let $x,x'\in\alpha$.  By Corollary \ref{cor:strong-relation-after-ev-period} there exists an $(x,x')$-walk in $R$ of length equal to $(p+1)p$. This means $x'\in (R^{p+1})^p(x)$. Hence, $[x']_{\sim_R}\bar{R}^p [x]_{\sim_R}$ and $\alpha\leftrightarrow_{\bar{R}}\alpha$, which proves that $\alpha\in \bar{X}_{\bar{R}}$. The right-to-left inclusion comes from the definition of recurrent set of $\bar{R}$. 

Notice that $\bar{R}$ restricted to a strongly connected component of $\bar{R}$ is a map. Indeed, suppose that there are $\alpha,\beta\in\bar{X}$, $\alpha\neq\beta$, such that $\alpha\in\bar{R}|_{[\gamma]_{\bar{R}}}(\gamma)$ and $\beta\in\bar{R}|_{[\gamma]_{\bar{R}}}(\gamma)$ for some $\gamma\in\bar{X}$. This means that $\gamma\bar{R}|_{[\gamma]_{\bar{R}}}\alpha$, $\gamma\bar{R}|_{[\gamma]_{\bar{R}}}\beta$ and for any $x\in\gamma, y\in\alpha, z\in\beta$ we have $x (R|_{\bigcup[\gamma]_{\bar{R}}})^{p+1} y$ and $x (R|_{\bigcup[\gamma]_{\bar{R}}})^{p+1} z$. Therefore, there are an $(x,y)$-walk and an $(x,z)$-walk of $R|_{\bigcup[\gamma]_{\bar{R}}}$, both of length equal to $p+1$. Hence, $y,z$ belong to the same class of $\sim_R$, $\alpha=\beta$, a contradiction.

Using a similar argument as in the paragraph above one can prove that $\bar{R}$ restricted to a strongly connected component is injective. In order to show that $\bar{R}$ restricted to a strongly connected component is surjective let $\alpha\in\bar{X}$ and take $\beta\in[\alpha]_{\bar{R}}$. Consider $\gamma = (\bar{R}|_{[\alpha]_{\bar{R}}})^{p-1}(\beta)$, where $p$ is an eventual period of $\bar{R}$ (see Lemma \ref{lem:R-bar-property}). We have $\bar{R}|_{[\alpha]_{\bar{R}}}(\gamma) = (\bar{R}|_{[\alpha]_{\bar{R}}})^p(\beta)$. By Proposition \ref{prop:identity-in-p-plus-q-power} and Corollary \ref{cor:strong-relation-after-ev-period} we get $\id_{[\alpha]_{\bar{R}}}\subset (\bar{R}|_{[\alpha]_{\bar{R}}})^p$. Therefore, $(\bar{R}|_{[\alpha]_{\bar{R}}})^p(\beta) = \beta$ and $\bar{R}|_{[\alpha]_{\bar{R}}}(\gamma) = \beta$. Thus, $\bar{R}|_{[\alpha]_{\bar{R}}}$ is a bijection.

Careful inspection of the proof of Lemma \ref{lem:R-bar-property} indicates that the variable $i$ may be replaced by any positive integer, which proves (iii).

Isomorphisms between objects $(X,R)$ and $(\bar{X},\bar{R})$ in the Szymczak category are given by Theorem \ref{thm:arisom}.
\qed

\begin{prop}
\label{prop:canonical-form-construction}
Let $(X,R)\in\Endo(\FRel)$ be in canonical form. Then $(X_R/_{\sim_R},\bar{R})$ is also in canonical form, where $\bar{R}$ is given as in (\ref{eq:R-bar}). Moreover, $(X,R)$ and $(X_R/_{\sim_R},\bar{R})$ are conjugate objects of $\Endo(\FRel)$. 
\end{prop}
\proof
By Theorem \ref{thm:each-obj-has-its-canonical-form}, the object $(X_R/_{\sim_R},\bar{R})$ is in canonical form. 

Consider the map $f\colon X\to X_R/_{\sim_R}$ such that $f(x):=[x]_{\sim_R}$. Since $X=X_R$, map $f$ is well-defined. Notice that for each $x\in X$ we have $\card [x]_{\sim_R} = 1$. Indeed, suppose on contrary that there are $x,x'\in[x]_{\sim_R}$ such that $x\neq x'$. Then there are an $(x,x)$-walk and an $(x,x')$-walk . This means that for some $y$ lying on both walks $\card R|_{[x]_R}(y)>1$, but $R|_{[x]_R}$ is a bijection, a contradiction. 

Using the above fact one can easily prove that $f$ is a bijection. By Proposition \ref{prop:FRel-auto}, map $f$ is an isomorphism between $X$ and $X_R/_{\sim_R}$ in $\FRel$.

We will show that $f\circ R = \bar{R}\circ f$. Let $p\in\NN_1$ be an eventual period of $R$ and let $x\in X$. We have 
$$
\begin{array}{rcl}
\bar{R}(f(x)) &=& \bar{R}([x]_{\sim_R}) = \{[y]_{\sim_R}\ |\ x R^{p+1} y\} = f(\{y\ |\ x R^{p+1} y\}) \\
 &=& f(\{y\ |\ y \in R^{p+1}(x)\}) = f(R^{p+1}(x)) = f(R(x)), \\
\end{array}
$$ 
which proves that $(X,R)$ and $(X_R/_{\sim_R},\bar{R})$ are conjugate in $\Endo(\FRel)$.
\qed

Because of Proposition \ref{prop:canonical-form-construction}, an object $(X,R)$ in canonical form is also said to be {\em canonical}.

\begin{figure}[h]
\begin{center}
  \includegraphics[width=0.85\textwidth]{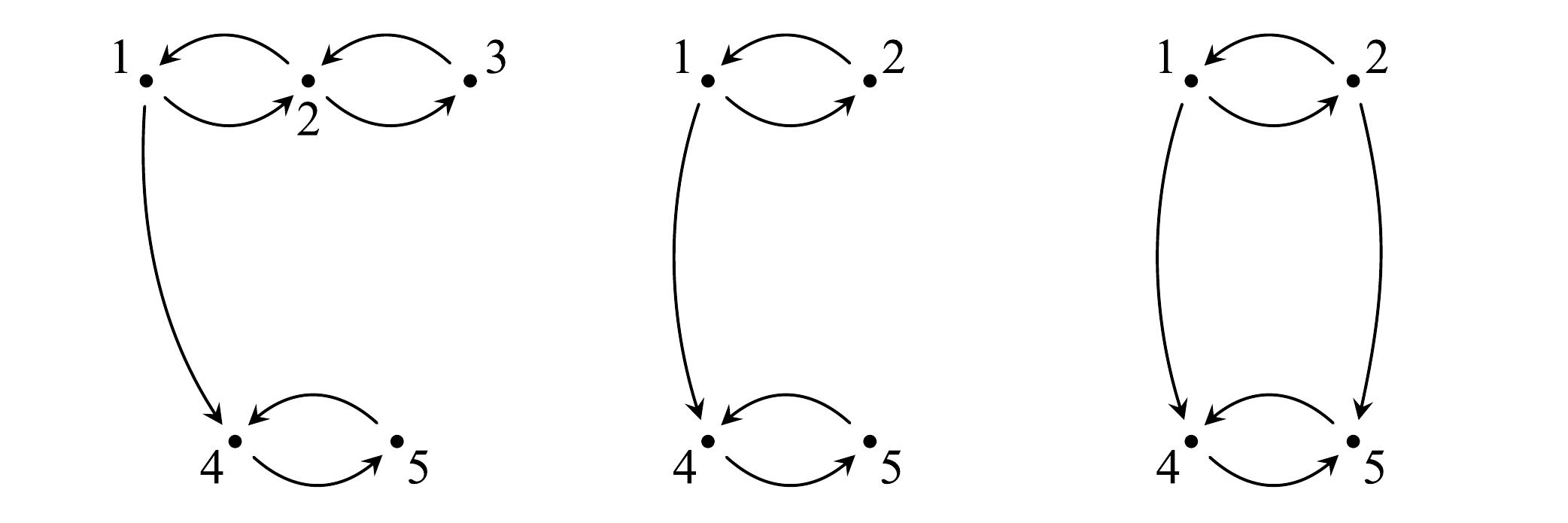}
\end{center}
  \caption{Relations $R_1$, $R_2$ and $R_3$ (from left to right) isomorphic in Szymczak's category. Only relation $R_3$ is in canonical form. }
  \label{fig:canon_mid}
\end{figure}

As an example we will show that the relation $R_1$ in Figure \ref{fig:canon_mid} is isomorphic to relation the $R_3$ in Szymczak's category. From now on, for the matrix representation A of a relation R we use the convention $A_{ij} = 1$ if $x_iRx_j$ and $0$ otherwise. We have 
$$
R_1 = 
\left(
\begin{array}{ccccc}
 0 & 1 & 0 & 1 & 0 \\
 1 & 0 & 1 & 0 & 0 \\
 0 & 1 & 0 & 0 & 0 \\
 0 & 0 & 0 & 0 & 1 \\
 0 & 0 & 0 & 1 & 0 \\
\end{array}
\right)
\qquad\text{ and }\qquad
R_3 = 
\left(
\begin{array}{cccc}
 0 & 1 & 1 & 0 \\
 1 & 0 & 0 & 1 \\
 0 & 0 & 0 & 1 \\
 0 & 0 & 1 & 0 \\
\end{array}
\right).
$$
An eventual period of $R_1$ is $p=4$. Relation $R_1$ has two strongly connected components $[1]_{R_1}:=\{1,2,3\}$ and $[4]_{R_1}:=\{4,5\}$, where the vertex number is also the row-column number of the matrix representation of the relation $R_1$. Moreover, we have $[1]_{\sim_{R_1}}=\{1,3\}$, $[2]_{\sim_{R_1}}=\{2\}$, $[4]_{\sim_{R_1}}=\{4\}$ and $[5]_{\sim_{R_1}}=\{5\}$. Using the formulas from the proof of Theorem \ref{thm:arisom} we get 
$$
T:=
\left(
\begin{array}{ccccc}
 1 & 0 & 1 & 0 & 1 \\
 0 & 1 & 0 & 1 & 0 \\
 0 & 0 & 0 & 1 & 0 \\
 0 & 0 & 0 & 0 & 1 \\
\end{array}
\right)
\qquad\text{ and }\qquad
S:=
\left(
\begin{array}{cccc}
 1 & 0 & 0 & 1 \\
 0 & 1 & 1 & 0 \\
 1 & 0 & 0 & 1 \\
 0 & 0 & 1 & 0 \\
 0 & 0 & 0 & 1 \\
\end{array}
\right).
$$
It is easy to check that $R_1\circ T = T\circ R_3$, $R_3\circ S = S\circ R_1$, $S\circ T = R_3^p$, and $T\circ S = R_1^p$. Therefore, $[S,p]$ and $[T,p]$ are mutually inverse isomorphisms and, in consequence, $(\{1,2,3,4,5\},R_1)$ and $(\{1,2,4,5\},R_3)$ are isomorphic in $\Szym(\FRel)$.

Let $(X,R)\in\Endo(\FRel)$ be an object in canonical form. The relation $R$ induces a partial order $\leq_R$ in $X/_{\leftrightarrow_R}$ defined by
\begin{equation}
\label{eq:induced-po}
[x]_R\leq_R [y]_R :\iff \text{there exists a } (y,x)\text{-walk in } R.
\end{equation}
Indeed, reflexivity of $\leq_R$ is obvious. If $[x]_R\leq_R [y]_R$ and $[y]_R\leq_R [x]_R$, then there are a $(y,x)$-walk and an $(x,y)$-walk. Hence, $x$ and $y$ are strongly connected, $[x]_R=[y]_R$. One can easily prove that $\leq_R$ is transitive. 

If $[x]_R\leq_R [y]_R$, then we say that the component $[y]_R$ is {\em connected with} the component $[x]_R$.

Now we present a few technical lemmas which give us information on isomorphisms in $\Szym(\FRel)$.
\begin{lem}
\label{lem:preserving-connections}
Let $(X,R), (X',R')\in\Endo(\FRel)$ be objects in canonical form isomorphic in $\Szym(\FRel)$. If $[S,k]\colon (X,R) \to (X',R')$ and $[T,k']\colon (X',R')\to (X,R)$ are mutually inverse isomorphisms, then for each $x\in X$ and each $y\in (T\circ S)(x)$ either $[y]_R$ is not connected with $[x]_R$ or $[x]_R=[y]_R$. 
\end{lem}
\proof
Let $x\in X$ and $y\in (T\circ S)(x)$. Since $[T\circ S,k'+k]=[\id_X,0]$, for some $l\in\NN$ we have $R^l(y)\subset R^l(T(S(x))) = R^{l+k+k'}(x)$. We can find $z\in [y]_R$ such that $z\in R^l(y)$. Therefore, $[x]_R$ is connected with $[y]_R$. This excludes that $[y]_R$ is connected with $[x]_R$, because $[x]_R=[y]_R$ would have to be the case.
\qed

\begin{lem}
\label{lem:preserving-components}
Let $(X,R), (X',R')\in\Endo(\FRel)$ be objects in canonical form isomorphic in $\Szym(\FRel)$. If $[S,\alpha]\colon (X,R) \to (X',R')$ is an isomorphism and $U$ is a component of $R$, then $S(U)$ contains uniquely determined component $V$ of $R'$ with the same period as $U$ such that no other component of $R'$ with non-empty intersection with $S(U)$ is connected with $V$.
\end{lem}
\proof
Let $[T,\beta]\colon (X',R')\to (X,R)$ be an isomorphism inverse to $[S,\alpha]$, let $x\in X$ and $U:=[x]_R$. Assume that the period of $R|_{[x]_R}$ is equal to $q$. Consider $S(x)$. Let $\cC_{S(x)}:=\{[x']_{R'}\ |\ S(x)\cap [x']_{R'}\neq\emptyset,\ x'\in X'\}$. 

We claim that there exists exactly one component $V\in\cC_{S(x)}$ such that no other component $W\in\cC_{S(x)}$ is connected with $V$. 

We have $S(x)\neq\emptyset$, because for some $t\in\NN$ we get $R^t(T(S(x)))=R^{t+\alpha+\beta}(x)$ and $R^{t+\alpha+\beta}(x)\cap [x]_R\neq \emptyset$. Hence, also $\cC_{S(x)}\neq\emptyset$. Since $\leq_{R'}$ is partial order, we take maximal elements of $\cC_{S(x)}$. Without loss of generality assume that there exist two components $[e]_{R'},[d]_{R'}\in\cC_{S(x)}$ which are maximal elements of $\leq_{R'}$ in $\cC_{S(x)}$ such that $e\in [e]_{R'}\cap S(x)$ and $d\in [d]_{R'}\cap S(x)$.

We have $T(e)\subset T(S(x))$ and $T(d)\subset T(S(x))$. Let us take components $[y]_R,[z]_R$ such that $y\in[y]_R\cap T(e)\neq\emptyset$ and $z\in[z]_R\cap T(d)\neq\emptyset$ and no other component of $R$ with non-empty intersection with $T(e)$ and $T(d)$, respectively, is connected with $[y]_R$ and $[z]_R$, respectively. We can take such components by an argument similar to the one in the paragraph above. 

Since there are no components in $R$ with non-empty intersection with $T(e)$ and $T(d)$ connected with $[y]_R$ and $[z]_R$, respectively, by Lemma \ref{lem:preserving-connections}, we get that $[y]_R=[x]_R$ and $[z]_R=[x]_R$.  

Since $y\in T(e)\cap [x]_R$ and there is $s\in\NN$ such that $x\in R^s(y)$, we have
$$
S(x)\subset S(R^s(y))\subset S(R^s(T(e))).
$$
Hence, $d\in S(x)\subset R'^s(S(T(e)))$. Moreover, there exists $t\in\NN$ such that
$$
R'^{t}(d)\subset R'^s(S(T(R'^{t}(e)))) = R'^{s+t+\alpha+\beta}(e).
$$
There is also an element $c\in R'^{t}(d)$ in $[d]_{R'}$. In consequence, $c\in R'^{s+t+\alpha+\beta}(e)$ and this means that $[d]_{R'}\leq_{R'} [e]_{R'}$, which contradicts the choice of components  $[e]_{R'}$ and $[d]_{R'}$. Therefore, $[e]_{R'} = [d]_{R'}$. Thus, we proved that there is only one component with the desired property.


Now, let $V$ be the component of $R'$ with the property that no other component of $R'$ from $\cC_{S(x)}$ is connected with $V$. Take $e\in V\cap S(x)$. We will prove that the component $V$ has the same period as $U$ (equal to $q$). Note that $x\in R^q(x)$, and then $e\in S(x)\subset S(R^q(x)) = R'^q(S(x))$. Hence, 
$$
e\in S(x)\subset R'^{q}(S(x))\subset R'^{2q}(S(x))\subset\ldots.
$$
Therefore, the period of $V$ is equal to either $k:=rq$ for some $r\in\NN_1$ or $k\in\NN_1$ and $k|q$.

As we proved above, $U$ is the component of $R$ with non-empty intersection with $T(e)$  such that no other component of $R$ with non-empty intersection with $T(e)$ is connected with it. Take $y\in T(e)\cap U$. Since we have the sequence of inclusions 
$$
y\in T(e)\subset R^{k}(T(e))\subset R^{2k}(T(e))\subset\ldots,
$$
the period of $U$ is equal to either $q = sk$ for some $s\in\NN_1$ or $q|k$. Combining the cases for the period of $V$ and $U$, we have to consider four cases. 

In the first case $q=srq$ it follows that $sr=1$ and $k=q$. In the second case $q|k$ and $k|q$, we also get immediately $k = q$. Consider the next case $q = sk$ and $k|q$. Since $e\in R'^{k}(e)$, we get $T(e)\subset R^{k}(T(e))$ and $y\in T(e)\cap U$. Therefore, either $y\in R^{k}(y)$ and then $k=q$ or there is $z\in T(e)\cap U$ such that $y\neq z$ and $y\in R^{k}(z)$. We have $y,z\in T(e)\subset R^{k}(T(S(x)))$. Also $T\circ S\circ R^l = R^{l+\alpha+\beta}$ for some $l\in\NN$. Hence, assuming without loss of generality that $l-k>0$ we get $R^{l-k}(y)\subset R^{l+\alpha+\beta}(x)$ and $R^{l-k}(z)\subset R^{l+\alpha+\beta}(x)$.

We have $x,y,z\in U$, $y\neq z$ and $R|_{U}$ is a bijection on $U$. There exist $y',z'\in U$ such that $y'\in R^{l-k}(y)$, $z'\in R^{l-k}(z)$ and $y'\neq z'$. Therefore, $y'\in R^{l+\alpha+\beta}(x)$ and $z'\in R^{l+\alpha+\beta}(x)$. That means $y'=z'$. This contradicts the choice of $y'$ and $z'$, so the alternative in the case cannot hold.

Analogously, it can be proved that in the fourth case $k=q$. Hence, the period of component $V$ is equal to $q$.


Now we will prove that $V\subset S(U)$. Let $e\in S(x)\cap V$, where $x\in U$ and $d\in R'(e)\cap V$, $y\in R(x)\cap U$. Suppose to the contrary that $d\notin S(y)$, that is there exist $w\in R(x)$, $w\in [w]_R\neq U$ such that $d\in S(w)$. Obviously, $U$ is connected with $[w]_R$. Since the period of $V$ is equal to $q$, $e\in R'^{q-1}(d)$ holds and $e\in R'^{q-1}(d)\subset R'^{q-1}(S(w))$. We have 
$$
T(e)\subset R^{q-1}(T(S(w))),
$$ 
and by repeating the reasoning of this proof we show that there exists $z\in T(e)\cap U$ such that  $z\in R^{q-1}(T(S(w)))$. Hence, $[w]_R$ is connected with $U$. By the assumption $U$ is connected with $[w]_R$. Therefore, $U=[w]_R$, a contradiction. Repeating the reasoning for each  element of $V$ we get $V\subset S(U)$. Since $V$ is uniquely determined by elements of  $U$ and no other component with non-empty intersection with $S(U)$ is connected with $V$, the proof is completed. 
\qed

\begin{lem}
\label{lem:isomorphisms-preserve-connecting-relation}
An isomorphism in $\Szym(\FRel)$ between objects $(X,R)$, $(X',R')\in\Endo(\FRel)$ in canonical form induces a bijection between $X/_{\leftrightarrow_R}$ and $X'/_{\leftrightarrow_{R'}}$. Moreover, the bijection maps $\leq_R$ to $\leq_{R'}$.
\end{lem}
\proof
First we prove that an isomorphism preserves connections between the corresponding components.

Let $[S,\alpha]\colon (X,R) \to (X',R')$, $[T,\beta]\colon (X',R')\to (X,R)$ be mutually inverse isomorphisms in $\Szym(\FRel)$. Let $U$ and $V$ be components of $R$ with periods $q_U$ and $q_V$, respectively. Let $W\subset S(U)$ and $Q\subset S(V)$ be components of $R'$ with periods $q_U$ and $q_V$, respectively, mentioned in Lemma \ref{lem:preserving-components}. Assume that $V\leq_R U$. We will prove that $Q\leq_{R'}W$. 

Take $e\in W$. There is an $x\in T(e)$ such that $x\in U$. Since $V\leq_R U$, there exists $y\in R^k(x)$ for some $k\in\NN_1$, where $y\in V$. We have $S(x)\subset S(T(e))$ and $S(R^k(x))\subset S(T(R'^k(e)))$. Hence, for some $l\in\NN$ we get
$$
R'^{l}(S(y))\subset R'^{l}(S(R^k(x)))\subset  R'^{k+l+\alpha+\beta}(e).
$$
Since $S(y)$ contains elements of $Q$ and no other component of $R'$ intersecting $S(y)$  is connected with $Q$, we take an element $d\in S(y)\cap Q$ and $c \in R'^{l}(d)\cap Q$. Therefore, $c\in R'^{k+l+\alpha+\beta}(e)$, that is $W$ is connected with $Q$, that is $Q\leq_{\bar{R}} W$.

Define a map $f\colon X/_{\leftrightarrow_R}\to X'/_{\leftrightarrow_{R'}}$ such that $f(U):=W$, where $W\subset S(U)$ and no other component of $R'$ intersecting $S(U)$  is connected with $W$. Since such a $W$ is determined uniquely (see Lemma \ref{lem:preserving-components}), the map $f$ is well-defined.

We will prove that $f$ is injective. Let $f(U)=W=f(V)$. 
Then $W\subset S(U)\cap S(V)$ and $T(W)\subset T(S(U))$, $T(W)\subset T(S(V))$. There is an $x\in T(W)\cap U$, where $U$ is the component of $R$ not connected with any other component intersecting $T(W)$. Similarly, $y\in T(W)\cap V$, where no other component intersecting $T(W)$ is connected with $V$. That means $U$ is connected with $V$ and $V$ is connected with $U$, hence $U=V$.

We prove that $f$ is surjective. Assume to the contrary that there is $W\in X'/{\leftrightarrow_{R'}}$ such that for each $U\in X/_{\leftrightarrow_R}$ the inequality $f(U)\neq W$ holds. We have $V\subset T(W)$ for some $V\in X/_{\leftrightarrow_R}$ and no other component of $R$ intersecting  $T(W)$ is connected with $V$. Since $S(V)\subset S(T(W))$, we get $W\subset S(V)$ and no other component intersecting with $S(V)$ is connected with $W$. Hence, $f(V)=W$, a contradiction. Therefore, the map $f$ is surjective.

In particular, $\card X/_{\leftrightarrow_R} = \card X'/_{\leftrightarrow_{R'}}$. By above facts we get that for each $U,V\in X/_{\leftrightarrow_R}$ if $U\leq_R V$, then $f(U)\leq_{R'} f(V)$. This proves that $f$ maps $\leq_R$ to $\leq_{R'}$.
\qed

\begin{cor}
\label{cor:same-number-of-components}
Isomorphic objects in $\Szym(\FRel)$ have the same number of components with the same periods. 
\end{cor}
\proof
Since for each object in $\Endo(\FRel)$ we can find an object in canonical form (see Theorem \ref{thm:each-obj-has-its-canonical-form}) isomorphic to the given one in $\Szym(\FRel)$, the composition of isomorphisms in $\Szym(\FRel)$ is an isomorphism between canonical forms. The conclusion comes from Lemmas \ref{lem:isomorphisms-preserve-connecting-relation} and \ref{lem:preserving-components}.
\qed

\begin{cor}
\label{cor:restricted-szym-isomorphism-is-bijective}
Let $(X,R),(X',R')\in\Endo(\FRel)$ be in canonical form and let $[S,\alpha]\colon (X,R)\to(X',R')$ be an isomorphism in $\Szym(\FRel)$. Assume that $f\colon X/_{\leftrightarrow_R}\to X'/_{\leftrightarrow_{R'}}$ is a bijection given by Lemma \ref{lem:isomorphisms-preserve-connecting-relation}. Then for each $x\in X$ the restriction $S|_{[x]_R\times f([x]_R)}$ is a bijection.
\end{cor}
\proof
By Lemma \ref{lem:preserving-components}, $f([x]_R)\subset S([x]_R)$ and the components $[x]_R$ and $f([x]_R)$ have the same periods. The relations $R$ and $R'$ restricted to these components respectively are bijections. Hence, $\card[x]_R = \card f([x]_R)$. 

We will prove that $S|_{[x]_R\times f([x]_R)}$ is a map. Let $[T,\beta]\colon (X',R')\to(X,R)$ be an inverse isomorphism to $[S,\alpha]$. Suppose that there exists $x\in X$ such that $\card S|_{[x]_R\times f([x]_R)}(x)>1$ and pick $y\in S|_{[x]_R\times f([x]_R)}(x)$. Then for each $x'\in[x]_R$ we have $\card S|_{[x]_R\times f([x]_R)}(x')>1$. Let $t\in T(y)$ such that $t\in [x]_R$. Then $S(t)\subset S(T(y))$ and $y',z'\in S(t)$, $y'\neq z'$ and $y',z'\in [y]_{R'}$. There are $y''\neq z''$ such that $y'',z''\in [y]_{R'}$ and
$$
y'',z''\in R'^{k}(S(t))\subset R'^k(S(T(y)))=R'^{k+\alpha+\beta}(y)
$$
for some $k\in\NN$. This yields $y''=z''$, a contradiction.

Since $\card[x]_R = \card f([x]_R)$, the map $S|_{[x]_R\times f([x]_R)}$ is a bijection.
\qed

\begin{lem}
\label{lem:restricted-isomorphism-commutativity}
Let $(X,R),(X',R')\in\Endo(\FRel)$ be in canonical form and let $[S,\alpha]\colon (X,R)\to(X',R')$ be an isomorphism in $\Szym(\FRel)$. Then for each $x\in X$
$$
R'\circ S|_{[x]_R\times f([x]_R)} = S|_{[x]_R\times f([x]_R)}\circ R.
$$
\end{lem}
\proof
Let $p\in\NN_1$ be an eventual period of $R$. Let us take $x'\in R'(S|_{[x]_R\times f([x]_R)}(x))$. By Corollary \ref{cor:restricted-szym-isomorphism-is-bijective}, there exist a $y'=S|_{[x]_R\times f([x]_R)}(x)$ and $x'\in R'(y')$. Consider $[x']_{R'}$. By Lemma \ref{lem:preserving-components}, there exists a $z\in X$ such that $[z]_R = f^{-1}([x']_{R'})$. Since $S|_{[z]_R\times f([z]_R)}$ is a bijection, assume that $x'=S|_{[z]_R\times f([z]_R)}(z)$. 

We will show that $z\in R(x)$. Notice that we have $S|_{[z]_R\times f([z]_R)}(z)\in R'(S|_{[x]_R\times f([x]_R)}(x))$ and
$$
T(S|_{[z]_R\times f([z]_R)}(z))\subset T(R'(S|_{[x]_R\times f([x]_R)}(x))) = R(T(S|_{[x]_R\times f([x]_R)}(x))).
$$
It follows that there is a $t\in T(S|_{[z]_R\times f([z]_R)}(z))$ such that $t\in [z]_R$ and $t\in R(T(S|_{[x]_R\times f([x]_R)}(x)))$. In particular, $t\in R(T(S(x)))$, therefore $R^k(t)\subset R(R^{\alpha+\beta+k}(x))$ for some $k\in\NN$. Let us take $\tilde{x}\in R^{\alpha+\beta+k}(x)$ such that $\tilde{x}\in [x]_R$. Since $R=R^{p+1}$, there exists a $\tilde{t}\in[z]_R$ such that $\tilde{t}\in R^k(t)$ and $\tilde{t}\in R(\tilde{x})$. 

Notice that $\bar{x}\in T(S|_{[x]_R\times f([x]_R)}(x))$ such that $\bar{x}\in[x]_R$ is uniquely determined by $x$, because the restrictions of $S$ and $T$ to the components are bijections. In consequence, $\tilde{x}=\bar{x}$. 

Furthermore, $t\in T(S|_{[z]_R\times f([z]_R)}(z))$ such that $t\in[z]_R$ is also uniquely determined by $z$. Let us take $\bar{t}\in R^{\alpha+\beta+k}(z)$ such that $\bar{t}\in[z]_R$. Then $\bar{t}\in R^k(t)$ and $\bar{t}=\tilde{t}$. 

Since $\tilde{t}\in R^k(t)$ and $\tilde{t}\in R^{\alpha+\beta+k}(z)$, we get $t\in R^{\alpha+\beta}(z)$. To sum up, we have $\tilde{x}\in R^{\alpha+\beta+k}(x)$, $\tilde{t}\in R(\tilde{x})$ and $z\in R^{mp-\alpha-\beta-k}(\tilde{t})$ for $mp>\alpha+\beta+k$ and $m\in\NN_1$. Combining these we get
$$
z\in R^{mp-\alpha-\beta-k}(R(R^{\alpha+\beta+k}(x))),
$$
which means that $z\in R^{mp+1}(x)$. Hence, $z\in R(x)$.

Since $x' = S|_{[z]_R\times f([z]_R)}(z)$ and $z\in R(x)$, we have $R'\circ S|_{[x]_R\times f([x]_R)} \subset S|_{[x]_R\times f([x]_R)}\circ R$. The proof of the opposite inclusion is analogous. 
\qed

\begin{lem}
\label{lem:restricted-rel-commutativity}
Let $(X,R)\in\Endo(\FRel)$ be in canonical form. Then for any $n\in\NN_1$ and for each $x\in X$
$$
R\circ R|_{[x]_R}^n = R|_{[x]_R}^n\circ R.
$$
\end{lem}
\proof
Let $y\in R(R|_{[x]_R}(x))$. There exists a $z=R|_{[x]_R}(x)$ such that $y\in R(z)$. There exists a $z'\in[y]_R$ such that $y\in R(z')$ and $y=R|_{[y]_R}(z')$. We have $y\in R^2(x)$ and $z'\in R^{p-1}(y)$, where $p$ is an eventual period of $R$. Thus, $z'\in R^{p-1}(y)\subset R^{p+1}(x) = R(x)$. Hence, $z'\in R(x)$ and $R\circ R|_{[x]_R} \subset R|_{[x]_R}\circ R$. The proof of the opposite inclusion is analogous.

Now assume that $R\circ R|_{[x]_R}^n = R|_{[x]_R}^n\circ R$. We have
$$
R\circ R|_{[x]_R}^{n+1} = R\circ R|_{[x]_R}^{n}\circ R|_{[x]_R} = R|_{[x]_R}^{n} \circ R\circ R|_{[x]_R} = R|_{[x]_R}^{n}\circ R|_{[x]_R}\circ R = R|_{[x]_R}^{n+1}\circ R.
$$
This completes the proof. 
\qed

\begin{thm}[Theorem \ref{thm:main2}]
\label{thm:szym-isom-induces-can-obj-endo-isom}
Let $(X,R),(X',R')\in\Endo(\FRel)$ be in canonical form. The objects $(X,R)$ and $(X',R')$ are isomorphic in $\Szym(\FRel)$ if and only if $(X,R)$ and $(X',R')$ are isomorphic in $\Endo(\FRel)$.
\end{thm}
\proof
Let $[S,\alpha]\colon (X,R)\to(X',R')$ and $[T,\beta]\colon (X',R')\to(X,R)$ be mutually inverse isomorphisms in $\Szym(\FRel)$ and let $t\in\NN_1$ be such that $T\circ S\circ R^t=R^{\alpha+\beta+t}$. Let us define morphisms $U\colon (X,R)\to(X',R')$ and $V\colon (X',R')\to(X,R)$ in $\Endo(\FRel)$ by
$$
\begin{array}{c}
U(x):= S|_{[x]_R\times f([x]_R)} (R|_{[x]_R}^{mp-\alpha-t}(x)), \\
V(x'):= T|_{[x']_{R'}\times f^{-1}([x']_{R'})} (R'|_{[x']_{R'}}^{mp-\beta}(x')), \\
\end{array}
$$
where $p\in\NN_1$ is an eventual period of $R$, $mp>\alpha+\beta+t$ for some $m\in\NN_1$, $f^{-1}$ is the inverse of bijection from Lemma \ref{lem:isomorphisms-preserve-connecting-relation} and $x\in X$, $x'\in X'$. We claim that $U$ and $V$ are mutually inverse isomorphisms in $\Endo(\FRel)$. 

By Corollary \ref{cor:restricted-szym-isomorphism-is-bijective}, both $U$ and $V$ are bijections. Using Lemma \ref{lem:restricted-isomorphism-commutativity} one can prove that $V(U(x))=R|_{[x]_R}^p(x)$. By Theorem \ref{thm:partition}, we have $R|_{[x]_R}^p(x)=[x]_{\sim_R}$ and $\card [x]_{\sim_R} = 1$ because $(X,R)$ is in canonical form. Therefore, $V(U(x))=R|_{[x]_R}^p(x)=\id_X(x)$. Similarly, one proves that $U(V(x'))=\id_{x'}(x')$.

Equalities $R'\circ U = U\circ R$ and $V\circ R' = R\circ V$ easily come from Lemmas \ref{lem:restricted-isomorphism-commutativity} and \ref{lem:restricted-rel-commutativity}.

Now let $F\colon (X,R)\to(X',R')$ and $G\colon(X',R')\to(X,R)$ be mutually inverse isomorphisms in $\Endo(\FRel)$. This means that $F\circ G = \id_{X'}$ and $G\circ F=\id_{X}$. Moreover, $F\circ R = R'\circ F$ and $G\circ R'=R\circ G$. Consider the morphisms $[F,p]\colon (X,R)\to(X',R')$ and $[G,p]\colon(X',R')\to(X,R)$ in $\Szym(\FRel)$. Using the facts above one can easily check that $[F\circ G,2p]=[\id_{X'},0]$ and $[G\circ F,2p]=[\id_{X},0]$.
\qed

%
\subsection{Classifying graphs}

Let $(X,R)\in\Endo(\FRel)$ be in canonical form. Define the map $l_{[x]_R}\colon [x]_R\times [x]_R\to\NN$ on connected components of $R$ such that
$$
l_{[x]_R}(x',x''):=m\mod q_{[x]_R}, \text{ if } x''\in R|_{[x]_R}^m(x'),
$$
where $q_{[x]_R}$ is the period of $R|_{[x]_R}$. Since the restriction $R|_{[x]_R}$ is a bijection and $(R|_{[x]_R})^k=(R|_{[x]_R})^{k+q_{[x]_R}}$ holds for $k\in\NN_1$, the maps $l_{[x]_R}$ are well-defined for each component $[x]_R$ of $R$.

Let $[x]_R$ and $[y]_R$ be  components of $R$ and let $q_{[x]_R}$ and $q_{[y]_R}$ be the periods of $R|_{[x]_R}$ and $R|_{[y]_R}$, respectively. Define the relation $\sim_{[x]_R[y]_R}\subset ([x]_R\times [y]_R)^2$ such that for $(x',y'),(x'',y'')\in [x]_R\times [y]_R$ we have
\begin{equation}
\label{eq:connections-between-components}
\begin{array}{rcl}
&(x',y')\sim_{[x]_R[y]_R} (x'',y'') :\iff & \\
& l_{[x]_R}(x',x'')=l_{[y]_R}(y',y'')\mod \gcd(q_{[x]_R},q_{[y]_R}).& \\
\end{array}
\end{equation}

\begin{prop}
\label{prop:connections-between-components-equivalence-relation}
The relation $\sim_{[\tilde{x}]_R[\tilde{y}]_R}$ on $[\tilde{x}]_R\times [\tilde{y}]_R$ is an equivalence relation for all components $[\tilde{x}]_R\neq [\tilde{y}]_R$ of $R$.
\end{prop}
\proof
Reflexivity of $\sim_{[\tilde{x}]_R[\tilde{y}]_R}$ is obvious. Let $(x,y)\sim_{[\tilde{x}]_R[\tilde{y}]_R} (x',y')$. Then $l_{[\tilde{x}]_R}(x',x)=q_{[\tilde{x}]_R}-l_{[\tilde{x}]_R}(x,x')$ and $l_{[\tilde{y}]_R}(y',y)=q_{[\tilde{y}]_R}-l_{[\tilde{y}]_R}(y,y')$. Since $q_{[\tilde{x}]_R}=q_{[\tilde{y}]_R}=0\mod \gcd(q_{[\tilde{x}]_R},q_{[\tilde{y}]_R})$ and $l_{[\tilde{x}]_R}(x,x')=l_{[\tilde{y}]_R}(y,y')\mod \gcd(q_{[\tilde{x}]_R},q_{[\tilde{y}]_R})$, we get $(x',y')\sim_{[\tilde{x}]_R[\tilde{y}]_R} (x,y)$. Hence, $\sim_{[\tilde{x}]_R[\tilde{y}]_R}$ is symmetric.

In order to prove transitivity of $\sim_{[\tilde{x}]_R[\tilde{y}]_R}$ let $(x,y)\sim_{[\tilde{x}]_R[\tilde{y}]_R} (x',y')$ and $(x',y')\sim_{[\tilde{x}]_R[\tilde{y}]_R} (x'',y'')$. Since $l_{[\tilde{x}]_R}(x,x')=l_{[\tilde{y}]_R}(y,y')\mod\gcd(q_{[\tilde{x}]_R},q_{[\tilde{y}]_R})$ and $l_{[\tilde{x}]_R}(x',x'')=l_{[\tilde{y}]_R}(y',y'')\mod\gcd(q_{[\tilde{x}]_R},q_{[\tilde{y}]_R})$, also 
$$
l_{[\tilde{x}]_R}(x,x')+l_{[\tilde{x}]_R}(x',x'')=l_{[\tilde{y}]_R}(y,y')+l_{[\tilde{y}]_R}(y',y'')\mod\gcd(q_{[\tilde{x}]_R},q_{[\tilde{y}]_R}). 
$$
It follows that $l_{[\tilde{x}]_R}(x,x'')=l_{[\tilde{y}]_R}(y,y'')\mod\gcd(q_{[\tilde{x}]_R},q_{[\tilde{y}]_R})$ and in consequence $(x,y)\sim_{[\tilde{x}]_R[\tilde{y}]_R} (x'',y'')$.
\qed

Note that $\sim_{[x]_R[y]_R}$ gives a partition of $[x]_R\times [y]_R$ into $\gcd(q_{[x]_R},q_{[y]_R})$ equivalence classes. Let $(X,R),(X',R')\in\Endo(\FRel)$ be in canonical form and $[S,\alpha]\colon(X,R)\to(X',R')$ be an isomorphism between the objects in $\Szym(\FRel)$. If $f([x]_R)\subset S([x]_R)$ and $f([y]_R)\subset S([y]_R)$ are components of $R'$ from Lemma \ref{lem:preserving-components}, where $f$ is the bijection from Lemma  \ref{lem:isomorphisms-preserve-connecting-relation}, then $\sim_{f([x]_R)f([y]_R)}$ on $f([x]_R)\times f([y]_R)$ also defines the partition into $\gcd(q_{[x]_R},q_{[y]_R})$ number of equivalence classes.

For $(X,R)$ in canonical form define the {\em number of connections between components $[x]_R$ and $[y]_R$ of $R$} as
\begin{equation}
\label{eq:number-of-connections}
\begin{array}{rcl}
&l_{[x]_R[y]_R}(R):=& \\
&\card\{ [(x',y')]_{\sim_{[x]_R[y]_R}}\in {[x]_R\times [y]_R/}_{\sim_{[x]_R[y]_R}}\ |\ (x',y')\in R|_{[x]_R\times [y]_R}\}. & \\
\end{array}
\end{equation}
This number determines how many equivalence classes of $\sim_{[x]_R[y]_R}$ are realized by connections between $[x]_R$ and $[y]_R$ components of $R$. The following proposition holds.

\begin{prop}
\label{prop:isomorphisms-preserve-number-of-connections}
Let $(X,R),(X',R')\in\Endo(\FRel)$ be in canonical form. If the objects are isomorphic in $\Szym(\FRel)$ and $f$ is the bijection between components of $R$ and $R'$ from Lemma \ref{lem:isomorphisms-preserve-connecting-relation}, then $l_{[x]_R[y]_R}(R) = l_{f([x]_R)f([y]_R)}(R')$.
\end{prop}
\proof
Let $[S,\alpha]\colon (X,R) \to (X',R')$ and $[T,\beta]\colon (X',R')\to (X,R)$ be mutually inverse isomorphisms. Consider components $[x]_R$ and $[y]_R$ and let $q_{[x]_R}$ and $q_{[y]_R}$ be the periods of $R|_{[x]_R}$ and $R|_{[y]_R}$, respectively. Let $\tilde{x}\in [x]_R$ and $e\in S(\tilde{x})\cap f([x]_R)$. Take all $e_1,\ldots,e_{k'}\in f([x]_R)$ such that $[(e,e_l)]_{\sim_{f([x]_R)f([y]_R)}}\neq [(e,e_m)]_{\sim_{f([x]_R)f([y]_R)}}$ for all $l\neq m$, $l,m\in\{1,\ldots,k'\}$. There exists a sequence $s'_1,\ldots,s'_{k'}\in\NN_1$ such that $e_l\in R'^{s'_l}(e)$ and $s'_l\neq s'_m\mod\gcd(q_{[x]_R},q_{[y]_R})$ for each $l\neq m$. In other words,  $l_{f([x]_R),f([y]_R)}(R') = k'$. 

We have also $T(e_l)\subset T(R'^{s'_l}(e))$. Take $x_l\in T(e_l)\cap [y]_R$ for each $l=1,\ldots,k'$. Then there is $t\in\NN$ such that for each $x_l$ we have $x_l\in R^{t+\alpha+\beta+s'_l}(\tilde{x})$. Since $s'_l\neq s'_m\mod\gcd(q_{[x]_R},q_{[y]_R})$ for $l\neq m$, we get $l_{[x]_R[y]_R}(R)\geq k'$. 

Assume to the contrary that there exist $x_1,x_2\in [y]_R$ such that the classes $[(\tilde{x},x_1)]_{\sim_{[x]_R[y]_R}}\neq [(\tilde{x},x_2)]_{\sim_{[x]_R[y]_R}}$ and for $e'\in S(x_1)$ and $e''\in S(x_2)$ we have $[(e,e')]_{\sim_{f([x]_R)f([y]_R)}} = [(e,e'')]_{\sim_{f([x]_R)f([y]_R)}}$. Then $e'\in R'^{s'}(e)$, $e''\in R'^{s''}(e)$ and $s'=s''\mod\gcd(q_{[x]_R},q_{[y]_R})$. Note that $e'\in R'^{s'}(S(\tilde{x}))$ and $e''\in R'^{s''}(S(\tilde{x}))$, hence $x_1\in T(e')\subset R^{t+\alpha+\beta+s'}(\tilde{x})$ and $x_2\in T(e'')\subset R^{t+\alpha+\beta+s''}(\tilde{x})$ for some $t\in\NN$. But $s'=s''\mod\gcd(q_{[x]_R},q_{[y]_R})$, so we get $[(\tilde{x},x_1)]_{\sim_{[x]_R[y]_R}} = [(\tilde{x},x_2)]_{\sim_{[x]_R[y]_R}}$, a contradiction. Therefore, $l_{[x]_R[y]_R}(R) = k' = l_{f([x]_R),f([y]_R)}(R')$.
\qed

Let $(X,R)\in\Endo(\FRel)$ and let $(\bar{X},\bar{R})\in\Endo(\FRel)$ be in canonical form such that the two objects are isomorphic in $\Szym(\FRel)$ (see Theorem \ref{thm:each-obj-has-its-canonical-form}). We define a {\em classifying graph $k(R)$}, that is a directed graph $k(R):= (V,E)$ such that $V:=\bar{X}/_{\leftrightarrow_{\bar{R}}}$ and $E:=\{([x]_{\bar{R}},[y]_{\bar{R}})\in V\times V\ |\ l_{[x]_{\bar{R}}[y]_{\bar{R}}}(\bar{R})\neq 0 \text{ and } [x]_{\bar{R}}\neq [y]_{\bar{R}}\}$. Vertices and edges of a classifying graph are labelled by positive integers. For an $[x]_{\bar{R}}\in V$ we label it by $\lab([x]_{\bar{R}}):= q_{[x]_{\bar{R}}}$, where $q_{[x]_{\bar{R}}}$ is the period of $\bar{R}|_{[x]_{\bar{R}}}$ and for an edge $([x]_{\bar{R}},[y]_{\bar{R}})\in E$ we label it by $\lab([x]_{\bar{R}},[y]_{\bar{R}}):=l_{[x]_{\bar{R}}[y]_{\bar{R}}}(\bar{R})$. Classifying graphs are invariants of isomorphic objects in $\Szym(\FRel)$. 

\begin{thm}
\label{thm:invariant}
Isomorphic objects in $\Szym(\FRel)$ have the same classifying graphs up to graph isomorphism preserving labels of vertices and edges. 
\end{thm}
\proof
By Theorem \ref{thm:arisom}, each object in $\Endo(\FRel)$ is isomorphic in $\Szym(\FRel)$ to some object in canonical form. Composing corresponding isomorphisms we get an isomorphism between canonical forms of the isomorphic objects. By Lemma \ref{lem:isomorphisms-preserve-connecting-relation}, Corollary \ref{cor:same-number-of-components} and Proposition \ref{prop:isomorphisms-preserve-number-of-connections} we get the proof. 
\qed

\begin{figure}[h]
\begin{center}
  \includegraphics[width=0.85\textwidth]{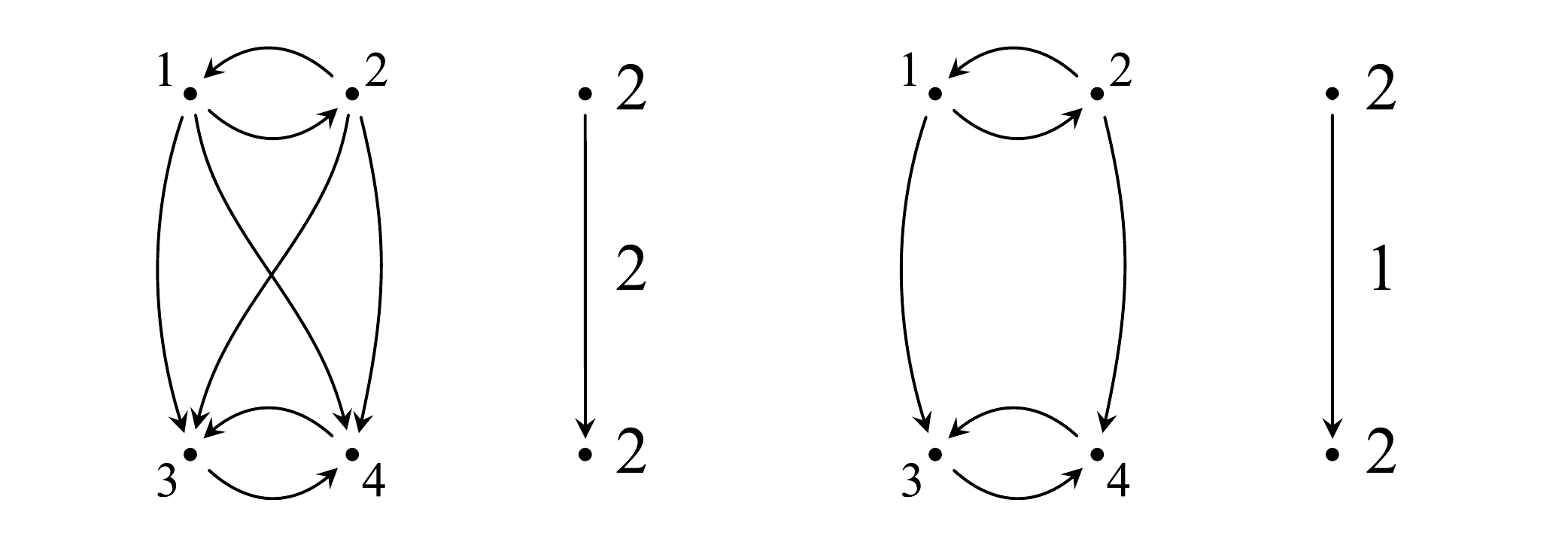}
\end{center}
  \caption{From left to right: relation $R$, classifying graph $k(R)$ of $R$, relation $R'$ and its classifying graph $k(R')$. The numbers of the vertices marked on relations digraphs denote the position in matrix representation of the relations. The numbers marked on the classifying graphs denote the labels of the vertices and the edges.}
  \label{fig:class_gr}
\end{figure}

Consider as an example objects $(X,R)$ and $(X',R')$ in canonical form given by
$$
R=
\left(
\begin{array}{cccc}
 0 & 1 & 1 & 1 \\
 1 & 0 & 1 & 1 \\
 0 & 0 & 0 & 1 \\
 0 & 0 & 1 & 0 \\
\end{array}
\right)
\qquad\text{ and }\qquad
R'=
\left(
\begin{array}{cccc}
 0 & 1 & 1 & 0 \\
 1 & 0 & 0 & 1 \\
 0 & 0 & 0 & 1 \\
 0 & 0 & 1 & 0 \\
\end{array}
\right).
$$
Both relations $R$ and $R'$ are pretty similar. They have two components with period both equal to $2$. One component is connected with the other. Assume that the first component in both relations is the set $\{1,2\}=:[1]_R=:[1]_{R'}$ and the second is the set $\{3,4\}=:[3]_R=:[3]_{R'}$ (numbers correspond to row-column positions of ones in matrix representation of these relations). We have $l_{[1]_R,[3]_R}(R)\neq 0$ and $l_{[1]_{R'},[3]_{R'}}(R')\neq 0$. More precisely, 
$$
\card ([1]_R\times [3]_R/_{\sim_{[1]_R,[3]_R}}) = \card ([1]_{R'}\times [3]_{R'}/_{\sim_{[1]_{R'},[3]_{R'}}}) = \gcd(2,2) = 2.
$$
By (\ref{eq:connections-between-components}), we easily compute that $l_{[1]_R,[3]_R}(R) = 2$ whereas  $l_{[1]_{R'},[3]_{R'}}(R') = 1$. By Theorem \ref{thm:invariant}, we conclude that $(X,R)$ and $(X',R')$ are not isomorphic in $\Szym(\FRel)$ (cf. Figure \ref{fig:class_gr}).

\begin{figure}[h]
\begin{center}
  \includegraphics[width=0.85\textwidth]{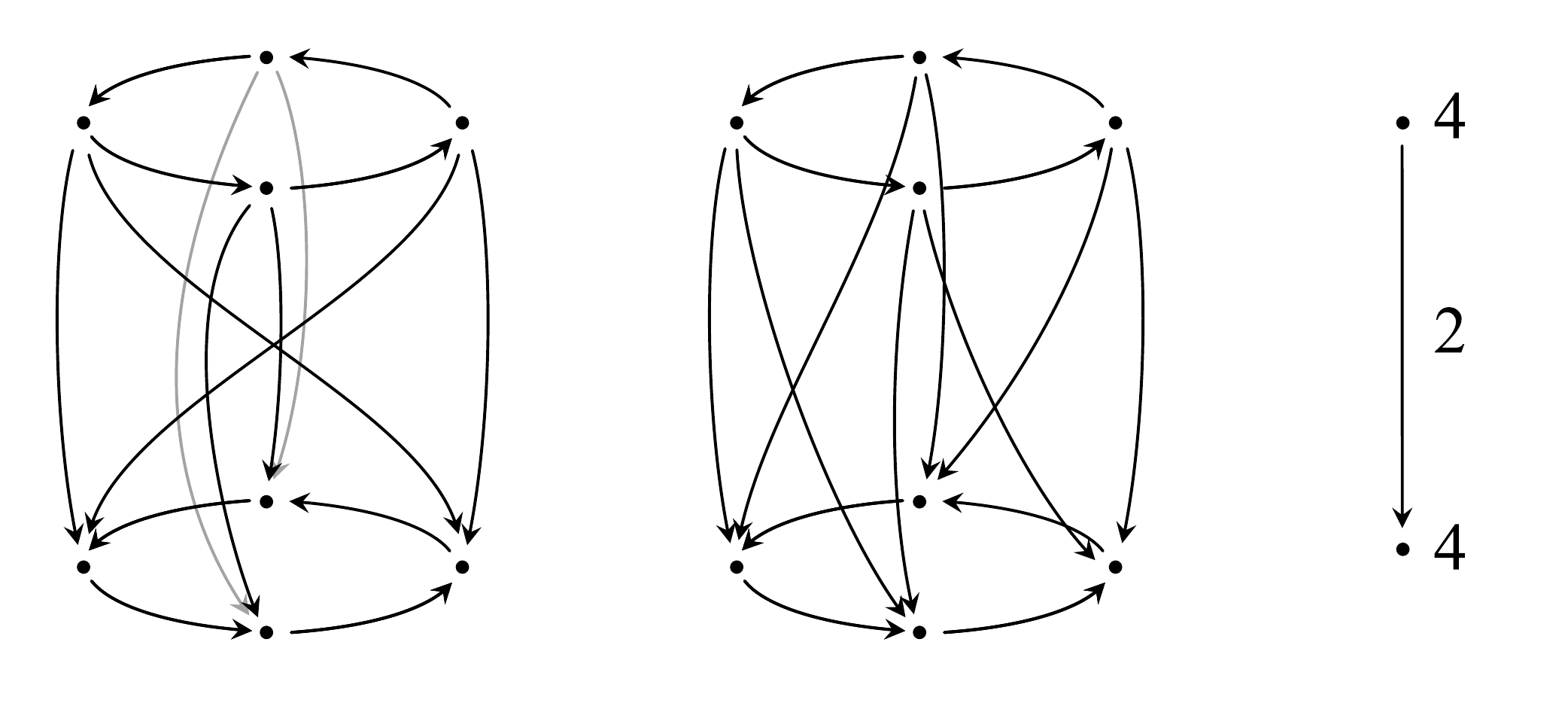}
\end{center}
  \caption{Two relations in canonical form with the same classifying graph (on the right) not isomorphic in $\Szym(\FRel)$}
  \label{fig:counterexample}
\end{figure}

Unfortunately, the classifying graph as an invariant of isomorphic object classes is not complete in sense that objects in $\Endo(\FRel)$ having the same classifying graphs up to graph isomorphism preserving labels of vertices and edges are isomorphic in $\Szym(\FRel)$. To see this, observe the example on Figure \ref{fig:counterexample}. Both relations are in canonical form and have the same classifying graphs but are neither isomorphic in $\Endo(\FRel)$ nor $\Szym(\FRel)$. 

%
\section{Final remarks}

The classification that we obtained allows us to distinguish non-isomorphic objects in $\Szym(\FRel)$ in an effective way. The main computational aspects involve strongly connected component detection, finding the period of a digraph component (the time complexity for both tasks is linear with respect to the sum of vertices and edges of the digraph; see \cite{JS1996}) and composition of relations (Boolean matrix multiplication). But in order to put this result into direct application in dynamics we need to consider relations with some algebraic structure, namely so-called linear relations. Recall, for vector spaces $X,Y$ over the field $\FF$ a relation $R\subset X\times Y$ is called {\em linear} (or {\em additive}; see \cite{ML1995}) if 
$$
\begin{array}{c}
(x_1,y_1)\in R, (x_2,y_2)\in R \implies (x_1+x_2,y_1+y_2)\in R, \\
(x_1,y_1)\in R \implies (ax_1,ay_1)\in R \text{ for each } a\in\FF. \\
\end{array}
$$
The sets with vector space structures with such relations constitute a category denoted by $\LRel$. 

We focus on linear relations since a multivalued generator of a dynamical system with non-acyclic values induces a linear relation. Such generators are common in sampled dynamics (see \cite{BMMP2020,HKMP2016}). Moreover, there are strong connections between $\LRel$ and $\FRel$. Therefore, we may use the $\Szym(\FRel)$ classification to understand $\Szym(\LRel)$. 

Notice that in general $\LRel$ is not a subcategory of the category of sets and relations since a given set may have more than one vector space structure. But there is a forgetful functor which forgets the linear structure of the space. Therefore, it is easy to check that if two objects equipped with relations on finite vector spaces are isomorphic in $\Szym(\LRel)$, then both objects are also isomorphic in $\Szym(\FRel)$. Thus, we may use the invariant from $\Szym(\FRel)$ as an invariant in $\Szym(\LRel)$. 

Consider the following example. Let $(\ZZ_3,R)$ and $(\ZZ_3,R')$ be objects of $\Endo(\LRel)$, where relations are defined in $\ZZ_3$ over $\ZZ_3$ with the standard operations. The relations are given by 
$$
R:=\{(0,0),(0,1),(0,2)\} \ \text{ and } \ R':=\{(0,0),(1,2),(2,1)\}.
$$ 
One can easily check that both relations are linear. Notice that relation $R$ is multivalued. After applying a functor induced by the forgetful functor we get two objects non-isomorphic in $\Szym(\FRel)$, because their classifying graphs are different (they have different numbers of components). Hence, $(\ZZ_3,R)$ and $(\ZZ_3,R')$ are non-isomorphic in $\Szym(\LRel)$. 

In such a way we may use the classification of $\Szym(\FRel)$ in understanding $\Szym(\LRel)$. On the other hand, the assumptions of a linear structure of relations is strong enough that it may significantly improve the classification of $\Szym(\LRel)$. For example, there are reasons to suppose that for linear relations over fields of finite (nonzero) characteristics the gradient structure of a relation between its components is no longer present or is trivial. Moreover, the stronger conditions imply that there are fewer morphisms in $\Szym(\LRel)$, so it is possible that the identification of two objects is not as common as in $\Szym(\FRel)$. Addressing these observations is beyond the scope of this paper and is a part of further research. We suppose that Szymczak's ideas may lead to the development of a Conley-index-type tool, enabling us to obtain dynamical information for systems reconstructed from data. 


\end{document}

%% file: MM_MathDefs.tex
\usepackage{amssymb}

\def\refeq#1{\if\workingver y(\ref{#1})-[[#1]]\else(\ref{#1})\fi}
\def\refth#1{\if\workingver y\ref{#1}-[[#1]]\else\ref{#1}\fi}
\def\mylabel#1{\if\workingver y\label{#1}{\bf\ \ [[#1]]\ \ }\else\label{#1}\fi}
\def\mybibitem#1{\if\workingver y\bibitem{#1}{\bf\ \ [[#1]]\ \
}\else\bibitem{#1}\fi}

\def\articletheorems{
\newtheorem{thm}{Theorem}[section]
\newtheorem{lem}[thm]{Lemma}

\newtheorem{cor}[thm]{Corollary}
\newtheorem{prop}[thm]{Proposition}

\newtheorem{algo}{Algorithm}[section] 
\newtheorem{alg}[thm]{Algorithm}  

}

\def\map{\rightarrow}

\newcommand{\mto}{\multimap}
\newcommand{\pto}{\nrightarrow}

\renewcommand{\emptyset}{\varnothing}
\renewcommand{\rho}{\varrho}
\renewcommand{\phi}{\varphi}
\renewcommand{\epsilon}{\varepsilon}

\def\cA{\text{$\mathcal A$}}
\def\cB{\text{$\mathcal B$}}
\def\cC{\text{$\mathcal C$}}

\def\cE{\text{$\mathcal E$}}

\newcommand{\id}{\operatorname{id}}

\newcommand{\dom}{\operatorname{dom}}

\newcommand{\card}{\operatorname{card}}

\newcommand{\im}{\operatorname{im}}

\renewcommand{\emptyset}{\varnothing}

\newcommand{\gim}{\operatorname{gim}}

\newcommand{\Inv}{\operatorname{Inv}}

\newcommand{\Per}{\operatorname{Per}}
\newcommand{\Endo}{\operatorname{Endo}}
\newcommand{\Auto}{\operatorname{Auto}}

\def\proof{{\bf Proof:\ }}

\def\begeq#1{\begin{equation}\mylabel{#1}}
\def\endeq{\end{equation}}

\def\mathobj#1{\mbox{$#1$}}

\def\FF{\mathobj{\mathbb{F}}}

\def\NN{\mathobj{\mathbb{N}}}

\def\ZZ{\mathobj{\mathbb{Z}}}

\def\implies{\;\Rightarrow\;}
\def\iff{\;\Leftrightarrow\;}

\def\setof#1{\mbox{$\{\,#1\,\}$}}

\def\0#1{\hbox{\kern25pt}$ #1 $\\}
\def\1#1{\hbox{\kern40pt}$ #1 $\\}
\def\2#1{\hbox{\kern55pt}$ #1 $\\}
\def\3#1{\hbox{\kern70pt}$ #1 $\\}

\newcounter{li}

\def\begalg#1{\begin{algo}\mylabel{#1}\normalshape:\small\baselineskip 10pt\\}
\def\endalg{\end{algo}}

\def\Figures(include=#1,cat=#2){
  \renewcommand{\textfraction}{.20}
  \renewcommand{\topfraction}{.80}
  \renewcommand{\bottomfraction}{.80}
  \renewcommand{\floatpagefraction}{.80}
  \newcount\figcount
  \figcount=0
  \let\includefigures=#1
  \def\figcat{#2}
}

\def\FigureFromFile[#1][#2](#3)#4
{
  \begin{figure}[htbp]
     \global\advance\figcount by 1
     \if\includefigures y\special{anisoscale #1.wmf, \the\hsize #2}\fi
     \vspace{#2}
     \caption{#4}
     \mylabel{#3}
   \end{figure}
}

\def\FigureFromFileTwoD[#1][#2,#3](#4)#5
{
  \begin{figure}[htbp]
     \global\advance\figcount by 1
     \if\includefigures y\special{anisoscale #1.wmf, #2 #3}\fi
     \vspace{#2}
     \caption{#5}
     \mylabel{#4}
   \end{figure}
}

\def\FigureF<#1>[#2](#3)#4
{
  \begin{figure}[htbp]
     \global\advance\figcount by 1
     \if\includefigures y\special{anisoscale \figcat/fig\number\figcount.wmf,
       \the\hsize #2}
     \fi
     \if\includefigures p
       \leavevmode
       \epsfxsize=\hsize
       \epsffile{#1}
     \fi
     \if\includefigures y
          \vspace{#2}
     \fi
     \caption{#4}
     \mylabel{#3}
   \end{figure}
}

\def\Figure[#1](#2)#3
{
  \begin{figure}[htbp]
     \global\advance\figcount by 1
     \if\includefigures y\special{anisoscale \figcat/fig\number\figcount.wmf,
       \the\hsize #1}
     \fi
     \if\includefigures p
       \leavevmode
       \epsfxsize=\hsize
       \epsffile{fig\number\figcount.eps}
     \fi
     \if\includefigures y
          \vspace{#1}
     \fi
     \caption{#3}
     \mylabel{#2}
   \end{figure}
}

%% file: MM_AlgDefs.tex
\def\0{\hbox{\kern5pt}}
\def\1{\hbox{\kern20pt}}
\def\2{\hbox{\kern35pt}}
\def\3{\hbox{\kern50pt}}
\def\4{\hbox{\kern65pt}}
\def\5{\hbox{\kern80pt}}
\def\6{\hbox{\kern95pt}}

